\documentclass[12pt]{article}
\usepackage{amsmath}
\usepackage{amssymb}
\usepackage{latexsym}
\usepackage[all]{xy}

\newtheorem{claim}{Claim}

\newtheorem{cor}{Corollary}

\newtheorem{lemma}{Lemma}
\newtheorem{prop}{Proposition}

\newtheorem{thm}{Theorem}

\newcounter{example} 
\newcounter{rem}
\newcommand{\Ann}{{\rm Ann}}
\newcommand{\ann}{{\rm ann}}

\newcommand{\barbbQ}{\ensuremath{\bar{\bbQ}}}
\newcommand{\barbbQp}{\ensuremath{\bar{\bbQ}_p}}

\newcommand{\bB}{\ensuremath{{\bf B}}}

\newcommand{\bbC}{\ensuremath{\mathbb{C}}}

\newcommand{\bbQ}{\ensuremath{\mathbb{Q}}}

\newcommand{\bbZ}{\ensuremath{\mathbb{Z}}}

\newcommand{\beq}{\begin{equation}}
\newcommand{\beql}[1]{\begin{equation}\label{#1}}
\newcommand{\bPf}{\noindent \textsc{Proof\ }}

\newcommand{\bT}{\ensuremath{{\bf T}}}

\newcommand{\cA}{\ensuremath{\mathcal{A}}}

\newcommand{\cF}{\ensuremath{\mathcal{F}}}

\newcommand{\charL}{\ensuremath{{\rm char}_\Lambda}}

\newcommand{\Cl}{{\rm Cl}}
\newcommand{\cL}{\ensuremath{\mathcal{L}}}

\newcommand{\cO}{\ensuremath{\mathcal{O}}}
\newcommand{\coker}{\ensuremath{{\rm coker}}}
\newcommand{\cores}{\ensuremath{{\rm cores}}}

\newcommand{\cR}{\ensuremath{\mathcal{R}}}
\newcommand{\cS}{\ensuremath{\mathcal{S}}}

\newcommand{\cV}{\ensuremath{\mathcal{V}}}

\newcommand{\displaymapdef}[5]
{\[
\begin{array}{rcrcl}
 #1 &:& #2 &\longrightarrow& #3 \\
    & &    &                    \\
    & & #4 &\longmapsto    & #5
\end{array}
\]}

\newcommand{\eeq}{\end{equation}}
\newcommand{\eg}{\emph{e.g.}}
\newcommand{\eitheror}[4] 
{\left\{                  
\begin{array}{ll}         
#1& \mbox{#2}\\
#3& \mbox{#4}
\end{array}
\right.}
\newcommand{\ePf}{\hspace*{\fill}~$\Box$\vertsp\par}
\newcommand{\eps}{\ensuremath{\varepsilon}}
\newcommand{\etc}{\emph{etc.}}

\newcommand{\fa}{\ensuremath{\mathfrak{a}}}
\newcommand{\fA}{\ensuremath{\mathfrak{A}}}

\newcommand{\fb}{\ensuremath{\mathfrak{b}}}
\newcommand{\fc}{\ensuremath{\mathfrak{c}}}
\newcommand{\fC}{\ensuremath{\mathfrak{C}}}
\newcommand{\fD}{\ensuremath{\mathfrak{D}}}
\newcommand{\fd}{\ensuremath{\mathfrak{d}}}

\newcommand{\fj}{\ensuremath{\mathfrak{j}}}
\newcommand{\fJ}{\ensuremath{\mathfrak{J}}}

\newcommand{\fm}{\ensuremath{\mathfrak{m}}}

\newcommand{\fP}{\ensuremath{\mathfrak{P}}}
\newcommand{\fq}{\ensuremath{\mathfrak{q}}}
\newcommand{\fQ}{\ensuremath{\mathfrak{Q}}}
\newcommand{\fr}{\ensuremath{\mathfrak{r}}}
\newcommand{\fR}{\ensuremath{\mathfrak{R}}}
\newcommand{\fs}{\ensuremath{\mathfrak{s}}}
\newcommand{\fS}{\ensuremath{\mathfrak{S}}}

\newcommand{\fw}{\ensuremath{\mathfrak{w}}}

\newcommand{\fX}{\ensuremath{\mathfrak{X}}}

\newcommand{\Gal}{{\rm Gal}}

\newcommand{\half}{\frac{1}{2}}
\newcommand{\Hom}{{\rm Hom}}
\newcommand{\ibid}{\emph{ibid.}}
\newcommand{\ie}{\emph{i.e.}}
\newcommand{\ignore}[1]{}
\newcommand{\im}{\ensuremath{{\rm im}}}

\newcommand{\inv}{^{-1}}

\newcommand{\ndiv}{\nmid}
\newcommand{\nin}{\not\in}
\newcommand{\ord}{{\rm ord}}

\newcommand{\rem}{\refstepcounter{rem}\noindent{\sc Remark \therem}}

\newcommand{\resp}{\emph{resp.}}

\newcommand{\tg}{\tilde{g}}
\newcommand{\Tr}{{\rm Tr}}
\newcommand{\vertsp}{\vspace{1ex}}

\newcommand{\jinf}{\ensuremath{\mathfrak{\fj_\infty}}}

\newcommand{\LG}{\ensuremath{\Lambda_\Gamma}}
\newcommand{\pnpo}{\ensuremath{{p^{n+1}}}}
\newcommand{\tC}{\ensuremath{\tilde{C}}}
\newcommand{\tCn}{\ensuremath{\tilde{C}_n}}
\newcommand{\tE}{\ensuremath{\tilde{E}}}
\newcommand{\tEn}{\ensuremath{\tilde{E}_n}}
\newcommand{\Vbar}[1]{\ensuremath{\bar{\cV}_{#1}}}
\newcommand{\Xinf}{\ensuremath{\mathfrak{\fX_\infty}}}

\newcommand{\Zpten}{\ensuremath{\bbZ_p\otimes}}
\begin{document}
\title{Some New Maps and Ideals in Classical Iwasawa Theory with Applications}
\author{D. Solomon
\\King's College, London\\
david.solomon@kcl.ac.uk}
\maketitle
\begin{abstract}
\noindent We introduce a new ideal $\fD$ of the $p$-adic Galois group-ring associated to a real
abelian field and a related ideal $\fJ$ for imaginary abelian fields. Both result from an equivariant, Kummer-type pairing applied to
Stark units in a $\bbZ_p$-tower of abelian fields and $\fJ$ is linked by
explicit reciprocity to a third ideal $\fS$ studied more generally in~\cite{laser}.
This leads to a new and unifying framework for the Iwasawa Theory of such fields
including a real analogue of Stickelberger's Theorem, links with certain Fitting ideals
and $\Lambda$-torsion submodules, and a new exact sequence related to the Main Conjecture.
\end{abstract}
\section{Introduction}
Let $K/k$ be any Galois extension of number fields and $p$ any {\em odd} prime number.
For each $n \geq -1$, we set $K_n=K(\mu_{\pnpo})$ and $G_n=\Gal(K_n/k)$. Let
$K_\infty=\bigcup_{n\geq -1}K_n$, $G_\infty=\Gal(K_\infty/k)$ and $\fX_\infty=\Gal(M_\infty/K_\infty)$  where
$M_\infty$ is the maximal abelian pro-$p$ extension of $K_\infty$ unramified outside
$p$. For each $n\geq 0$ we shall write $\Gamma_n$ for $\Gal(K_\infty/K_n)$.

Now suppose that $k=\bbQ$ and $\Gal(K/\bbQ)$ is abelian. Let
$K_n^+$ and $K_\infty^+$ be the maximal real subfields and let $G_n^+=\Gal(K_n^+/\bbQ)$. By applying
Kummer theory to a particular sequence of cyclotomic `units'
$\{\eps_n\}_{n\geq 0}$ of the fields $K_n^+$ we obtain a pair of related $G_\infty$-(semi)linear maps denoted
$\fd_\infty$ and $\fj_\infty$, defined on $\fX_\infty$ and taking values in the `plus' and `minus' parts respectively
of the completed group ring $\bbZ_p[[G_\infty]]$. The main purpose of this paper is to explore systematically the properties of
these maps and their images, the latter being ideals of $\bbZ_p[[G_\infty]]$ which we denote $\fD_\infty$ and $\fJ_\infty$ respectively.
In so doing, we establish precise links with, and/or applications to, the following areas, among others:
the Galois structure of the class group and also (units)/(cyclotomic units) for a real, absolutely abelian field;
Greenberg's and Vandiver's Conjectures; explicit reciprocity laws and the map `$\fs$' introduced in~\cite{twizo,laser};
the `$\Lambda$-torsion' submodule of  $\fX_\infty$; the Main Conjecture over $\bbQ$. (These
connections are dealt with in successive sections whose content is described in more detail below.)
As well as producing new results,
we therefore feel that these maps and ideals provide a
unifying approach to several significant
areas of the Iwasawa Theory of abelian number fields, and one that has so far been largely overlooked.
(We also mention in passing a technical advantage of this approach as compared to some others,
namely the way it works naturally at the the group-ring level.
This means that
we never need to decompose using $p$-adic characters of $\Gal(K/\bbQ)$ in the present paper.
Consequently no exceptions or special
treatments are necessary for the `non-semi-simple' case, \ie\ when $p$ divides $[K:\bbQ]$.)

In Section~\ref{sec: gencon} of this paper we consider general $K/k$ as above and define the basic Kummer-theoretic pairing between $\fX_\infty$ and
norm-coherent sequences of (global) $p$-units, taking values in $\bbZ_p[[G_\infty]]$.

In Section~\ref{sec: cunits and annihilation} and from there on, $K$ is almost always taken to be an (absolutely) abelian field and
$k=\bbQ$. The pairing applied to the above sequence $\{\eps_n\}_{n\geq 0}$ then produces both the map $\fd_\infty$,
and the map $\fj_\infty$ (which is its `mirror-twist') and hence the ideals
$\fD_\infty$ and $\fJ_\infty$. The images of the latter in $\bbZ_p[G_n]$ are ideals denoted $\fD_n$ and $\fJ_n$ respectively.
We give a concrete description of $\fD_n$ (modulo $\pnpo$)
using power-residue symbols. This dovetails perfectly with Thaine's methods to show that $\fD_n$ annihilates
the $p$-part of the class group $\Cl(K_n^+)$ (at least if $K_\infty/K_n$ is totally ramified). In this sense it
can be seen as an analogue in the plus-part of $\bbZ_p[G_n]$ of the ($p$-adified) Stickelberger ideal in the minus-part.
We also give an abstract characterisation of $\fD_n$ in terms of the $\bbZ_p[G_n]$-dual of $p$-units.

In Section~\ref{sec: var K and kernels} we study the behaviour of $\fj_\infty$ and $\fd_\infty$ when we replace $K$ by a subfield. This contributes
to the proof of the main result of this section, namely that the common kernel of $\fd_\infty$ and $\fj_\infty$ is precisely the subgroup
$\Gal(M_\infty/N^0_\infty)$ of $\fX_\infty$. Here, $N^0_\infty$ denotes the field obtained by adjoining to $K_\infty$ all $p$-power roots of those
units whose local absolute norms above $p$ are trivial. This implies in particular that if $K_\infty^+$ has only one prime above $p$ then
$\fd_\infty$ and $\fj_\infty$ are injective iff Greenberg's
conjecture holds for the extension $K_\infty^+/K^+$.

Section~\ref{sec: d and D for K=Q} gives more precise results in the case $K=\bbQ$ \ie\ $K_n=\bbQ(\mu_{\pnpo})$.
We show that $\bbZ_p[G_n^+]/\fD_n$ is then naturally isomorphic to the {\em Pontrjagin dual} of the $p$-part of the quotient of units by
cyclotomic units in $K_n^+$.
In particular, $\fD_n$ is precisely the Fitting ideal of this dual and also, by a result of Cornacchia and Greither,
that of the $p$-part of $\Cl(K_n^+)$. Further links between $\fd_\infty$ and the Conjectures of Greenberg and Vandiver then follow naturally.

Back in the case of general abelian $K$, Section~\ref{sec: inertia subgps} starts by
considering the restriction of $\fj_\infty$ to the product of inertia subgroups in $\fX_\infty$.
An explicit reciprocity law due to Coleman connects this restriction with the projective limit $\fs_\infty$ of
certain maps $\fs_n$ defined in terms of $p$-adic logarithms and complex $L$-values at $s=1$ for odd Dirichlet characters.
(The latter maps were introduced and studied in  a more general context in~\cite{twizo,laser}.)
This connection has many consequences. For instance, we show that $\ker(\fj_\infty)=\ker(\fd_\infty)=\Gal(M_\infty/N_\infty^0)$ is also precisely
the torsion submodule of $\fX_\infty$ as a module over the Iwasawa algebra $\Lambda$. We also deduce a new
4-term exact sequence of torsion $\Lambda$-modules, involving both $\fD_\infty$
and $\fS_\infty$ (the image of $\fs_\infty$).

For Section~\ref{sec: analysis of exact sequence} we return to the special case
$K=\bbQ$ of Section~\ref{sec: d and D for K=Q} and use results of~\cite{Iwasawa on some mods...} on the image of the $p$-adic logarithm to show that
$\fS_\infty$ is then precisely the limit of the ($p$-adified) Stickelberger ideals.
The above-mentioned  exact sequence then shows that for a given even,
non-trivial power $\omega^j$ of the Teichm\"uller character $\omega$, the Main Conjecture over $\bbQ$ can be rephrased as an equality
between the characteristic power series in $\Lambda$ of the $\omega^j$-components of $\ker(\fd_\infty)$ and $\coker(\fd_\infty)$.

We have become aware of a few previous papers containing constructions having something is common with our
$\fD_\infty$, $\fd_\infty$ and $\fj_\infty$. Of these, the closest to ours in spirit seems to be~\cite{Kraft-Schoof}.
Its aims are, however, much narrower than ours and relate principally to computing the structure of certain Iwasawa modules in
the case where $K$ is real-quadratic (and $p=3$, assumed not
to split in $K$. See also~\cite{Schoof} for different but related techniques and computations.)
In Remark~\ref{rem: correspondence with K-S} we outline the connections between some of the
results of~\cite{Kraft-Schoof} and some of those
in our Sections~\ref{sec: cunits and annihilation} and~\ref{sec: d and D for K=Q}.
Next, if we restrict to $K=\bbQ$ and let $m$ be an odd integer, then
the value of the so-called {\em Soul\'{e} character}
$\chi_m:\fX_\infty\rightarrow\bbZ_p$ at $h\in\fX_\infty^-$ turns
out to be simply the integral of the
$(1-m)$th power of the cyclotomic character with respect to $\fd_\infty(h)$, regarded as a $\bbZ_p$-valued measure on $G_\infty$.
See~\cite{I-S} p.~54. It might therefore be interesting to compare the results of
our Section~\ref{sec: analysis of exact sequence} with some of those mentioned in~\cite[\S 3]{I-S}.
Finally, we mention that, in a different context and for different purposes, Section~6.2 of \cite{Sharifi} contains the construction of
a map `$\phi_2$' which is related to our $\fj_\infty$ for cyclotomic fields $K$.

Looking to the future, the first problem is to generalise to any abelian
$K$ the results obtained in Sections~\ref{sec: d and D for K=Q} and~\ref{sec: analysis of exact sequence} for $K=\bbQ$.
Beyond this, one might like to consider abelian extensions $K/k$ with $r:=[k:\bbQ]>1$. For imaginary quadratic $k$,
elliptic units might substitute for the $\eps_n$'s but, in some ways,
a stronger analogy with the present case can be expected when $k$ is totally real and the $K_n$'s are CM.
The best available substitutes for the $\eps_n$'s are then `Rubin-Stark elements' for $K_n^+/k$.
Unfortunately, not only do these lie {\em a priori}
in a certain $r$th exterior power of $S$-units of $K_n^+$ tensored with $\bbQ$ but
their existence is only conjectural. It is, however, strongly supported by computations \eg~\cite[\S 3.5]{vercors} which also
suggest that for $n\geq 0$ one can use $p$-units and tensor only with $\bbZ_{(p)}=\bbQ\cap\bbZ_p$.
By assuming this, one could mimic some of Sections~\ref{sec: gencon}
to~\ref{sec: d and D for K=Q}, replacing $\fX_\infty$ by an appropriate
$r$th exterior power, \etc\ On the other hand, the map $\fs_n$ is already defined unconditionally in this case in~\cite{twizo, laser}.
To connect it with a generalised $\fj_\infty$, the Congruence Conjecture (formulated in~\cite{laser} and tested numerically in~\cite{vercors})
would be precisely the required substitute for Coleman's reciprocity law mentioned above.

My thanks are due to Jean-Robert Belliard, Ralph Greenberg and Romyar Sharifi for their comments and particularly
to Thong Nguyen Quang Do for helpful
discussions and for pointing out~\cite{Kraft-Schoof},~\cite{I-S} and the possibility of some connection with~\cite{NN}, yet to be determined.
I am also very grateful to the anonymous referee for his/her useful comments on an earlier version of this paper.

Some notation: If $F$ is any field then $\mu(F)$ denotes the group of roots of unity in $F^\times$ with the subgroup
$\mu_m(F)$ (\resp\ $\mu_{p^\infty}(F)$) consisting of those of order dividing $m>0$ (\resp\ of $p$-power order).
We write $\xi_m$ for the generator $\exp(2\pi i/m)$ of $\mu_m:=\mu_m(\bbC)$ and $\mu_{p^\infty}$ for $\mu_{p^\infty}(\bbC)$.
A `number field' $L$ is always a finite extension of $\bbQ$
contained it its algebraic closure $\barbbQ\subset\bbC$. We write $\cO_L$, $E(L)=\cO_L^\times$ and $S_r(L)$ respectively for its ring of integers, its
unit group and the set of its places (or prime ideals)
dividing an integer $r>1$. If $F/F'$ is an abelian extension of number fields and $\fq$ is a prime ideal of $\cO_{F'}$ unramified in $F$ then
$\sigma_{\fq,F/F'}$ denotes the (unique) Frobenius element attached to
$\fq$ in $\Gal(F/F')$.

\section{A General Construction}\label{sec: gencon}
Let $K/k$ be any Galois extension of number fields and let $K_n$ and other notations be as above. (Thus $K=K_{n}$ for $n=-1$
and possibly for some $n\geq 0$.)
Throughout this paper we shall write $\pi^m_n$ for the natural
restriction map $G_m\rightarrow G_n$ (where $m\geq n\geq -1$)
or indeed for the homomorphism of group rings $\cR [G_m]\rightarrow \cR[G_n]$ obtained by
$\cR$-linear extension, for any commutative ring $\cR$.
We identify $G_\infty$ with the projective limit of the $G_n$'s w.r.t. the $\pi^m_n$'s.

Let $P_n$ denote the maximal abelian extension of $K_n$ of exponent dividing $\pnpo$ and unramified
outside $S_p(K_n)$. (It contains $K_{2n+1}$ and is finite over $K_n$). Let
$M_n$ be the maximal abelian pro-$p$-extension of $K_n$ unramified outside $S_p(K_n)$. Thus $M_n$ contains $P_n K_\infty$ and
$M_\infty=\bigcup_{n\geq -1}M_n=\bigcup_{n\geq -1}P_n$. Both $M_\infty/k$ and $M_n/k$ (for all $n$) are (infinite) Galois extensions.
We write $\fX_n$ for the profinite group $\Gal(M_n/K_n)$ (so that $\fX_\infty={\displaystyle\lim_\leftarrow}\fX_n$)
and $\bar{\fX}_n$ for the quotient $\fX_n/\fX_n^\pnpo$ which identifies with $\Gal(P_n/K_n)$ and hence is finite.

Set $\cV_n=E_{S_p}(K_n):=\cO_{K_n,S_p(K_n)}^\times$ (the group of `$p$-units' of $K_n$)
and $\Vbar{n}:=\cV_n/\cV_n^\pnpo$.
Kummer theory gives a unique, well-defined pairing $\langle\ ,\ \rangle_n:\Vbar{n}\times \bar{\fX}_n\rightarrow \bbZ/\pnpo\bbZ$  satisfying
\[
h(\alpha^{1/\pnpo})/\alpha^{1/\pnpo}=\zeta_n^{\langle\bar\alpha,\bar{h}\rangle_n}\ \ \ \ \mbox{ for all $\alpha\in\cV_n$ and $h\in\fX_n$}
\]
where $\zeta_n$ denotes $\xi_\pnpo$ and $\alpha^{1/\pnpo}$ is any of the $\pnpo$th roots of $\alpha$ (all lying in $P_n$).
We abbreviate $\bbZ/\pnpo\bbZ$ to $\cR_n$ so that
$\langle\ ,\ \rangle_n$ is $\cR_n$-bilinear.
\begin{figure}[t]
\[
\xymatrix{
&&&\ \ \ \ \ M_\infty\ \ \ \ \ \ \ \ar@{-}[ld]\ar@/_1pc/@{-}[lldd]_{\fX_\infty}\\
&&M_n\ar@{-}[ld]\ar@/^1pc/@{-}[lddd]^{\fX_n}\\
&K_\infty\ar@{-}[dd]\ar@/_/@{-}[dd]_{\Gamma_n}\ar@/_2pc/@{-}[ldddd]_{G_\infty}\\
&& P_n\ar@{-}[ld]\ar@/^1pc/@{-}[ld]^{\bar{\fX}_n}\ar@{-}[uu]\\
&K_n\ar@{-}[d]\ar@/_/@{-}[ldd]_{G_n}\\
&\ \ \ \ \ \ \ \ K_{-1}=K\ar@{-}[ld]\ar@/^1pc/@{-}[ld]^{G_{-1}}\\
k\ar@{-}[d]\\
\bbQ
}
\]
\end{figure}
We shall write $\chi_{cyc}:G_\infty\rightarrow\bbZ_p^\times$ for  the
$p$-cyclotomic character, determined by $g(\zeta)=\zeta^{\chi_{cyc}(g)}$ for any
$g\in G_\infty$ and $\zeta\in\mu_{p^\infty}$, so that $\chi_{cyc}(\Gamma_n)\subset 1+\pnpo\bbZ_p$ for all $n\geq 0$.
Reducing $\chi_{cyc}$ modulo $\pnpo$ gives a character $\chi_{cyc,n}: G_n\rightarrow\cR_n^\times$ for all $n\geq 1$.
For any $g\in G_n$
and $h\in\fX_n$ we define $g.h$ to be $\tilde{g}h\tilde{g}\inv$ for any lift $\tilde{g}$ of $g$ to $\Gal(M_n/k)$.
This determines a left $G_n$-action on $\fX_n$, hence on $\bar{\fX}_n$, and it follows easily from the
definition that
\beq\label{eq:Gn-action}
\langle g\bar\alpha,g.\bar{h}\rangle_n=\chi_{cyc,n}(g)\langle \bar\alpha,\bar{h}\rangle_n
\ \ \ \ \mbox{for all $\alpha\in\cV_n$, $h\in\fX_n$ and $g\in G_n$.}
\eeq
Next, $\langle\ ,\ \rangle_n$ gives rise to a group-ring-valued pairing
$\{\ ,\ \}_n:\Vbar{n}\times \bar{\fX}_n\rightarrow \cR_n[G_n]$ defined by
\beq\label{eq: def of equivariant pairing}
\{\bar\alpha,\bar{h}\}_n=\sum_{g\in G_n}\langle \bar\alpha,g\inv.\bar{h}\rangle_ng=
\sum_{g\in G_n}\chi_{cyc,n}(g)\inv\langle g\bar\alpha,\bar{h}\rangle_n g
\eeq
for all $\alpha\in\cV_n$ and $h\in\fX_n$, which is $\cR_n[G_n]$-linear in the second variable and  $\cR_n[G_n]$-semilinear in the first.
More precisely, there is an involutive automorphism $\iota_n$ of
$\cR_n[G_n]$ sending $\sum_{g\in G_n}a_gg$ to $\sum_{g\in G_n}a_g\chi_{cyc,n}(g)g\inv$ and
equations~(\ref{eq:Gn-action}) and~(\ref{eq: def of equivariant pairing}) show that
\beq\label{eq: Gn-action for equiv pairing}
\{x\bar\alpha,y.\bar{h}\}_n=\iota_n(x) y\{\bar\alpha,\bar{h}\}_n\ \ \ \
\mbox{for all $\alpha\in\cV_n$, $h\in\fX_n$ and $x,y\in\cR_n[G_n]$.}
\eeq
Clearly, $\iota_n(\{\bar\alpha,\bar{h}\}_n)$
is $\cR_n[G_n]$-linear in $\bar\alpha$ and $\iota_n$-semilinear in $\bar{h}$ and~(\ref{eq: def of equivariant pairing}) gives
\beq\label{eq: formula for iotan applied to pairing}
\iota_n(\{\bar\alpha,\bar{h}\}_n)=\sum_{g\in G_n}\langle g\inv\bar\alpha,\bar{h}\rangle_n g
\eeq
If $m\geq n\geq -1$ then $M_m\supset M_n$ and we write $\rho^m_n$ for the
restriction $\fX_m\rightarrow\fX_n$ and $\bar{\rho}^m_n:\bar{\fX}_m\rightarrow \bar{\fX}_n$. We also write
$N^m_n$ for the norm map $K_m^\times\rightarrow K_n^\times$ inducing
$\bar{N}^m_n:\Vbar{m}\rightarrow \Vbar{n}$, and
$\bar{\pi}^m_n$ for the ring homomorphism
$\cR_m [G_m]\rightarrow\cR_n[G_n]$ which acts as
$\pi^m_n$ on the elements of $G_m$ and as the reduction $\cR_m\rightarrow\cR_n$ on the coefficients. From~(\ref{eq:Gn-action})
and the fact that $\chi_{cyc,m}(g)\equiv 1\pmod{\pnpo}$ for all $g\in\Gal(K_m/K_n)$, we deduce:
\begin{prop}\label{prop:norm compat for equiv pairing} If $m\geq n\geq -1$ then the diagram
\[
\xymatrix{
\Vbar{m}\times\bar{\fX}_m\ar[dd]_{\bar{N}^m_n\times \bar{\rho}^m_n}\ar[rrr]^{\{\ ,\ \}_m}&&&\cR_m[G_m]\ar[dd]^{\bar{\pi}^m_n}\\
&&&&&\\
\Vbar{n}\times\bar{\fX}_n\ar[rrr]^{\{\ ,\ \}_n}&&&\cR_n[G_n]\\
}
\]
commutes. \ePf
\end{prop}
Passing to projective limits with respect to $\bar{N}^m_n$, $\bar{\rho}^m_n$ and $\bar{\pi}^m_n$ for $m\geq n\geq 0$,
we obtain a pairing
\[
\{\ ,\}_\infty:=\lim_{\longleftarrow\atop n\geq 0}\{\ ,\}_n:\lim_{\longleftarrow\atop n\geq 0}\Vbar{n} \times \lim_{\longleftarrow\atop n\geq 0}\bar{\fX}_n\longrightarrow
\lim_{\longleftarrow\atop n\geq 0}\cR_n[G_n]
\]
Each of the last three limits above has another interpretation. The third identifies (as a compact topological ring) with the completed group-ring
$\Lambda_G:=\bbZ_p[[G_\infty]]={\displaystyle\lim_{\longleftarrow}}\ \bbZ_p[G_n]$.
For future reference, it may help to make this identification explicit:
Decompose $\bar{\pi}^m_n$ as $\beta_{m,n;n}\circ \phi^m_{n}$ where
$\phi^j_{k}: \cR_j[G_j]\rightarrow\cR_j[G_k]$ (for $j\geq k\geq 0$) is $R_j$-linear, acting as
$\pi^j_k$ on the elements of $G_j$, and $\beta_{i,j;k}: \cR_i[G_k]\rightarrow\cR_j[G_k]$ (for $i\geq j\geq k\geq 0$) simply reduces coefficients
modulo $p^{j+1}$. Then a sequence
$(x_n)_n$ of ${\displaystyle\lim_{\longleftarrow}}\ \cR_n[G_n]$ gives rise to a sequence
$y_k:=(\phi^j_{k}(x_j))_j$ for each $k\geq 0$ lying in the limit
${\displaystyle\lim_{\longleftarrow\atop}}\ \cR_j[G_k]$ w.r.t.\ the $\beta_{i,j;k}$'s ($k$ fixed). Thus we get a sequence
\[
(y_k)_k\in \lim_{\longleftarrow\atop k\geq 0}(\lim_{\longleftarrow\atop j\geq k}\ \cR_j[G_k])=\lim_{\longleftarrow\atop k\geq 0}(\bbZ_p[G_k])=\Lambda_G
\]
Conversely, an element $(z_k)_k\in{\displaystyle \lim_{\longleftarrow}}(\bbZ_p[G_k])$ gives rise to the element $(z_n \pmod{\pnpo})_n$ of
${\displaystyle\lim_{\longleftarrow}}\ \cR_n[G_n]$.
Similar decompositions of $\bar{\rho}^m_n$ and $\bar{N}^m_n$ respectively identify ${\displaystyle\lim_{\longleftarrow}}\ \bar{\fX}_n$ with
\[
{\displaystyle\lim_{\longleftarrow\atop k\geq 0}}({\displaystyle\lim_{\longleftarrow\atop j\geq k}}\ \fX_k/p^{j+1}\fX_k)=
{\displaystyle\lim_{\longleftarrow\atop k\geq 0}}\ \fX_k=\fX_\infty
\]
and ${\displaystyle\lim_{\longleftarrow}}\ \Vbar{n}$ with
\[
{\displaystyle\lim_{\longleftarrow\atop k\geq 0}}({\displaystyle\lim_{\longleftarrow\atop j\geq k}}\ \cV_k/p^{j+1}\cV_k)=
{\displaystyle\lim_{\longleftarrow\atop k\geq 0}}\ (\bbZ_p\otimes\cV_k)=:\cV_\infty
\]
Thus we may regard $\{\ ,\ \}_\infty$ as a continuous $\Lambda_G$-valued pairing between the compact, topological $\Lambda_G$-modules
$\cV_\infty$ and $\fX_\infty$. If
$\iota_\infty={\displaystyle\lim_{\longleftarrow}}\ \iota_n$ denotes the continuous, involutive automorphism of $\Lambda_G$
sending  $g\in G_\infty$ to $\chi_{cyc}(g)g\inv$, then~(\ref{eq: Gn-action for equiv pairing}) leads to
\beq\label{eq: Lambda-action for equiv pairing at infty}
\{x\underline{\alpha},y.h\}_\infty=\iota_\infty(x) y\{\underline{\alpha},h\}_\infty\ \ \ \
\mbox{for all $\underline{\alpha}\in\cV_\infty$, $h\in\fX_\infty$ and $x,y\in\Lambda_G$}
\eeq

Fix $n\geq -1$ and let $\fq$ be a prime of $K_n$ not dividing $p$. Then $\mu_\pnpo$ injects into $(\cO_{K_n}/\fq)^\times$, so that $\pnpo$ divides
$N\fq-1$ and
the image of $\mu_\pnpo$ is precisely the subgroup of $((N\fq-1)/\pnpo)$th powers in $(\cO_{K_n}/\fq)^\times$.
If $\beta\in K^\times$ is
a local unit at $\fq$ 
we write $\left\{\frac{\beta}{\fq}\right\}_n$ for the additive $\pnpo$th
power-residue symbol mod $\fq$, \ie\ the unique element of $\cR_n$ satisfying
\[
\beta^{(N\fq-1)/\pnpo}\equiv\zeta_n^{\left\{\frac{\beta}{\fq}\right\}_n}\pmod{\fq}
\]
Now any $\bar{h}\in\bar{\fX}_n$ can be written as
$\sigma_{\fq,P_n/K_n}$ for some such ideal $\fq$, so the following characterises the pairing
$\{\ ,\ \}_n$.
\begin{prop}\label{prop: pairing and power residue}
Let $n\geq-1$, let $\fq$ be a prime of $K_n$ not dividing $p$ and $\alpha\in\cV_n$. Then
\[
\{\bar{\alpha},\sigma_{\fq,P_n/K_n} \}_n=
\sum_{g\in G_n}\left\{\frac{\alpha}{g\inv(\fq)}\right\}_n g
\ \ \ \mbox{and}\ \ \
\iota_n(\{\bar{\alpha},\sigma_{\fq,P_n/K_n}\}_n)=\sum_{g\in G_n}\left\{\frac{g\inv(\alpha)}{\fq}\right\}_n g
\]
\end{prop}
\bPf\ A well-known argument gives
$\langle\bar{\alpha},\sigma_{\fq,P_n/K_n} \rangle_n=\left\{\frac{\alpha}{\fq}\right\}_n$ so the second equation follows
from~(\ref{eq: formula for iotan applied to pairing}).
Since also $g\inv.\sigma_{\fq,P_n/K_n}=\sigma_{g\inv(\fq),P_n/K_n}$ for all $g\in G_n$, the first follows
from~(\ref{eq: def of equivariant pairing}).
\ePf

\section{Cyclotomic Units and the Annihilation of Real Classes}\label{sec: cunits and annihilation}
We shall suppose henceforth that $k=\bbQ$ and $\Gal(K/k)$ is abelian so  that $K_n$ is an (absolutely) abelian field
for all $n\geq -1$. {\em We shall also suppose $n\in\bbZ$, $n\geq 0$} unless explicitly stated otherwise.
We write $c$ for the element of $\Gal(\barbbQ/\bbQ)$ induced by complex conjugation and also for its restriction
to $K_n$ for any $n$. In the notation of the Introduction, its fixed field is $K_n^+$, $G_n^+\cong G_n/\{1,c\}$
and $K_\infty^+=\bigcup_{n\geq -1}K_n^+$.
If $M$ is any module for one of the (commutative) rings
$\cR_n[G_n]$, $\bbZ_p[G_n]$, $\Lambda_G$, \etc\ we shall also write
$M^+$ (\emph{resp.} $M^-$) for the submodule of $M$ on which $c$ acts trivially  (\emph{resp.} by $-1$). Since $p\neq 2$, we
have $M=M^+ \oplus M^-$, corresponding to the decomposition
$m=m^++m^-:=\half(1+c)m+ \half(1-c)m$ for each $m\in M$.

For any abelian field $F$ we shall write $f_F$ for its conductor (\ie\ the smallest integer $f\geq 1$ such that $F\subset\bbQ(\mu_f)$) so a prime
number $r$ divides $f_F$ iff it ramifies in $F$.
If $F\neq \bbQ$  then $f\geq 3$ and we write
$\eps_F$ for the cyclotomic `unit' attached to $F$ namely $N_{\bbQ(\mu_{f_F})/F}(1-\xi_{f_F})\in \cO_F$.
We shall need the following result, see \eg\ Lemma~2.1 of~\cite{SolGalRel}.
\begin{lemma}\label{lemma: norm relations} Suppose $F,F'$ are abelian fields with $F\supset F'\supsetneq \bbQ$. Using an {\em additive}
notation for Galois action on  the mutiplicative group of $F'$
\[
N_{F/F'}\eps_F=\bigg(\prod_r(1-\sigma_{r,F'/\bbQ}\inv)\bigg)
\eps_{F'}
\]
where the product (in $\bbZ[\Gal(F'/\bbQ)]$) runs over prime the numbers $r$ dividing $f_F$ but not $f_{F'}$. Moreover if $f_F$ is a power of some prime number, say $r$, then
$N_{F/\bbQ}\eps_F=r$ (so $\eps_F$ is an $r$-unit of $F$). Otherwise $N_{F/\bbQ}\eps_F=1$ (so $\eps_F$ is a unit of $F$).\ePf
\end{lemma}
For brevity, we shall write $f_n$ for $f_{K_n}$ for any $n\geq -1$, so $f_n$ is the l.c.m.\ of $\pnpo$ and $f_{-1}$. If
$n\geq 0$ then $K_0\neq \bbQ$ and we set
\[
\eps_n:=N_{K_n/K_n^+}\eps_{K_n}=\eps_{K_n}^{1+c}\in \cV_n^+\ \ \ \mbox{and}\ \ \ \eta_n:=\half\otimes\eps_n\in(\bbZ_p\otimes\cV_n)^+
\]
(Notice that $\eps_n$ coincides with $\eps_{K_n^+}$
provided $f_n$ equals $f_{K_n^+}$. Since
$K_n=K_n^+(\mu_p)$, this  holds iff $p|f_{K_n^+}$, \eg\ if $p>3$ or $n>0$.) Let $\bar{\eta}_n$ denote the image of
$\eta_n$ in $(\bbZ_p\otimes\cV_n)/\pnpo(\bbZ_p\otimes\cV_n)$ which
identifies canonically with $\bar{\cV}_n$ so that $\bar{\eta}_n=\half\bar{\eps}_n\in \bar{\cV}_n^+$.
Taking $x=y=c$ in equation~(\ref{eq: Gn-action for equiv pairing}) gives
$\{\bar{\eta}_n,c.\bar{h}\}_n=-\{\bar{\eta}_n,\bar{h}\}_n$ for all $h\in\fX_n$ and hence
\beq\label{eq: pairing with eta depends only on minus part}
\{\bar{\eta}_n,\bar{h}\}_n=\{\bar{\eta}_n,\bar{h}^-\}_n\in \cR_n[G_n]^-
\eeq
so that
$\iota_n(\{\bar{\eta}_n,\bar{h}\}_n)\in \cR_n[G_n]^+$. The $\cR_n$-linear extension of the
restriction map $G_n\rightarrow G_n^+$ identifies the $\cR_n[G_n]^+$
with $\cR_n[G_n^+]$ so  $\Cl(K_n^+)/\pnpo\Cl(K_n^+)$ becomes an
$\cR_n[G_n]^+$-module.
\begin{thm}\label{thm: Stick in plus part mod pnpo}
If $n\geq 0$ and $h\in\fX_n$ then $\iota_n(\{\bar{\eta}_n,\bar{h}\}_n)$ annihilates $\Cl(K_n^+)/\pnpo\Cl(K_n^+)$.
\end{thm}
(We shall shortly use this to construct explicit annihilators without prior knowledge of $\fX_n$.)\vertsp\\
\bPf\ The Theorem will follow from the following, apparently much weaker statement.
\begin{claim}\label{claim: for thm 1}
Let $\fq$ be a prime  of $K_n$ dividing a rational prime $q$ which splits completely in $K_n$.
Then
the element $\half\sum_{g\in G_n}\left\{\frac{g\inv(\eps_n)}{\fq}\right\}_n g$ of $\cR_n [G_n]^+$
annihilates the image of the class $[\fq^+]$ in $\Cl(K_n^+)/\pnpo\Cl(K_n^+)$, where $\fq^+$ is the prime of $K_n^+$ below $\fq$.
\end{claim}
Assume this for the moment.
Let $H_n^+$ be the maximal unramified abelian extension of $K_n^+$ of exponent dividing $\pnpo$ so that
the Artin map defines an isomorphism $\Cl(K_n^+)/\pnpo\Cl(K_n^+)\rightarrow\Gal(H_n^+/K_n^+)$
sending the class of $\fc\in\Cl(K_n^+)$ to $\sigma_\fc$, say.
Since $p\neq 2=[K_n:K_n^+]$, the restriction map $\Gal(K_n.H_n^+/K_n)\rightarrow\Gal(H_n^+/K_n^+)$
is an isomorphism. Moreover, $K_n.H_n^+\subset P_n$ so the restriction map
$\phi:\bar{\fX}_n=\Gal(P_n/K_n)\rightarrow\Gal(H_n^+/K_n^+)$ factors through the previous one and is surjective.
But $K_n.H_n^+/K_n^+$  is abelian so
$c$ acts trivially on $\Gal(K_n.H_n^+/K_n)$
from which it follows that $\phi(\bar{\fX}_n^-)=\{0\}$  and $\phi(\bar{\fX}_n^+)=\Gal(H_n^+/K_n^+)$.
Now choose any $\fc\in\Cl(K_n^+)$ and any element $h\in\fX_n$. Since
$\bar{\fX}_n=\bar{\fX}_n^+\oplus\bar{\fX}_n^-$, \v{C}ebotarev's Theorem implies the existence of $\fq$ satisfying the hypotheses of
the Claim such that
$\sigma_{\fq,P_n/K_n}^-=\bar{h}^-\in\bar{\fX}_n^-$ {\em and} such that
$\sigma_{\fq,P_n/K_n}^+\in\bar{\fX}_n^+$ maps to $\sigma_\fc$ by $\phi$, hence so does
$\sigma_{\fq,P_n/K_n}$. On the one hand it follows from~(\ref{eq: pairing with eta depends only on minus part})
and~Proposition~\ref{prop: pairing and power residue} that
$\iota_n(\{\bar{\eta}_n,\bar{h}\}_n)=
\iota_n({\{\bar{\eta}_n,\sigma_{\fq,P_n/K_n}\}_n})=
\half\sum_{g\in G_n}\left\{\frac{g\inv(\eps_n)}{\fq}\right\}_n g$. On the other hand it follows from
the properties of the Artin map (and the fact that $\fq^+$ splits in $K_n$) that $\fc$ and  $[\fq^+]$
have the same image in $\Cl(K_n^+)/\pnpo\Cl(K_n^+)$. Thus the Claim
implies that $\iota_n(\{\bar{\eta}_n,\bar{h}\}_n)$ annihilates the image of $\fc$ for {\em all} $h$ and $\fc$.

Our proof of Claim~\ref{claim: for thm 1} is close to that of Thaine's Theorem as given in~\cite[\S 15.2]{Wash}.
The splitting condition implies $q\ndiv f_n$ and that $\mu_\pnpo$ injects into $(\cO_{K_n}/\fq)^\times=(\bbZ/q\bbZ)^\times$. In particular,
$\pnpo|(q-1)$ so we may choose a primitive root $t\in \bbZ$
modulo $q$ such that $t^{(q-1)/\pnpo}\equiv\zeta_n\pmod{\fq}$.
We denote by $K_{n,q}$ the field $K_n(\xi_q)$ which is easily seen to be unramified over $K_{n,q}^+$ at all finite primes. Since $q\neq 2$
both $K_{n,q}$ and $K_{n,q}^+$ have conductor $qf_n$. The extension $K_{n,q}/K_n$ is totally tamely ramified at all primes above $q$
hence so is $K_{n,q}^+/K_n^+$ (and $K_{n,q}/K_n$ is unramified elsewhere). Therefore
$\bbQ(\mu_{f_n})$ and $K_{n,q}^+$ are linearly
disjoint over $K_n^+$ with compositum $\bbQ(\mu_{qf_n})$.
$\Gal(K_{n,q}/K_n)$ identifies by restriction with $\Gal(K_{n,q}^+/K_n^+)$ and is cyclic of degree $q-1$
generated by $\tau:\xi_q\rightarrow\xi_q^t$.
Set $\eps_{n,q}=N_{\bbQ(\mu_{qf_n})/K_{n,q}^+}(1-\xi_{q}\xi_{f_n})$ which is clearly conjugate over $\bbQ$ to
$\eps_{K_{n,q}^+}$. Lemma~\ref{lemma: norm relations}
implies that it is a unit of $K_{n,q}^+$  and that $N_{K_{n,q}^+/K_n^+}(\eps_{n,q})=1$, since $\sigma_{q,K_n^+/\bbQ}=1$.
By Hilbert's Theorem~90,
we can therefore choose $\beta\in K_{n,q}^{+,\times}$
such that $\tau(\beta)/\beta=\eps_{n,q}$. It follows that
$\ord_{\fR^+}(\beta)=\ord_{\fR^+}(\tau^i(\beta))$ for any prime $\fR^+$ of $K_{n,q}^+$ and for $i=1,\ldots,q-1$,
so $\ord_{\fR^+}(N_{K_{n,q}^+/K_n^+}\beta)=(q-1)\ord_{\fR^+}(\beta)$. If $\fR^+\ndiv q$ then $\fR^+$ has
ramification index $1$ or $2$
over $K_n^+$ (the latter case if $K_n/K_n^+$ is ramified at the prime below $\fR^+$, which requires $\fR^+|p$ and $f_n$ a power of $p$).
We deduce that the principal fractional ideal of $K_n^+$ generated by
$N_{K_{n,q}^+/K_n^+}\beta$ is of the form $\fa I^{(q-1)/2}$ where $I$ is a fractional ideal prime to $q$ and $\fa$ has support above $q$.
For each $g\in G_n$, we write $\fQ_g$ for the unique prime of $K_{n,q}$ dividing $g(\fq)$ and $\fQ_g^+$ for the prime of $K_{n,q}^+$ below it
(dividing $g(\fq^+)$), so $\fQ_g$ is split over $\fQ^+_g$. We set
\[
a_g:=\ord_{g(\fq^+)}(\fa)=\frac{1}{q-1}\ord_{\fQ^+_g}(N_{K_{n,q}^+/K_n^+}\beta)=\ord_{\fQ^+_g}(\beta)=\ord_{\fQ_g}(\beta)\in\bbZ
\]
The stabiliser of $\fq^+$ in $G_n$ is $\{1,c\}$ so in $\Cl(K_n^+)$ we have
$\sum_{g\in G_n}a_g g [\fq^+]=2[\fa]=(1-q)[I]\in\pnpo\Cl(K_n^+)$. Thus the Claim, and hence the Theorem, will follow once we have proven that
\beq\label{eq: to prove Claim 1}
a_g\equiv\left\{\frac{g\inv(\eps_n)}{\fq}\right\}_n \pmod{\pnpo}\ \ \ \mbox{for every $g\in G_n$}
\eeq
Since $\ord_{\fQ_g}(1-\xi_q)=1$ we can write $\beta=(1-\xi_q)^{a_g}v$ for some $v\in K_{n,q}^\times$, a local unit at $\fQ_g$. Therefore
$\eps_{n,q}=((1-\xi_q^t)/(1-\xi_q))^{a_g}\tau(v)/v=(1+\xi_q+\ldots+\xi_q^{t-1})^{a_g}\tau(v)/v$ and since $\tau$ acts trivially on the residue field at
$\fQ_g$ by total ramification, we deduce $\eps_{n,q}\equiv t^{a_g}\pmod{\fQ_g}$. On the other hand, $1-\xi_q\xi_{f_n}$
is congruent to $1-\xi_{f_n}$ modulo all primes of $\bbQ(\mu_{qf_n})$ dividing $q$,
from which it follows easily that $\eps_{n,q}\equiv\eps_n\pmod{\fQ_g}$. Thus
$\eps_n\equiv t^{a_g}\pmod{g(\fq)}$ so that
$g\inv(\eps_n)^{(q-1)/\pnpo}\equiv t^{a_g(q-1)/\pnpo}\equiv \zeta_n^{a_g}\pmod{\fq}$, giving~(\ref{eq: to prove Claim 1}).\ePf
\rem\ One can in fact deduce Theorem~\ref{thm: Stick in plus part mod pnpo} from Theorem~1.3 of~\cite{Rubin-GUICG} (a
far more general elaboration of Thaine's method). This is explained briefly below.
A minor complication occurs if  $f_n$ is a power of $p$ but the main virtue of our {\em ab~initio} proof is its
much greater simplicity and directness compared to Rubin's proof of Theorem~1.3.
This is natural enough given the specialness  of our situation.

In Rubin's  Theorem~1.3, take
`$K$', `$F$', `$N$' and `$G$' in  to be $\bbQ$, $K_n^+$, $\pnpo$ and $G_n^+$ respectively. Let
his `$V$' and `$A$' be $E(K^+_n)/E(K^+_n)^\pnpo$ and $\Cl(K_n^+)/\pnpo\Cl(K_n^+)$ respectively. Given any
$\bar{h}\in\bar{\fX}_n$, we may take `$\alpha$' to be the map $\alpha_{\bar{h}}:v\mapsto \iota_n(\{v,\bar{h}\}_n)$.
It follows easily from~\cite[Lemma~1.6(ii)]{Rubin-GUICG} that Rubin's `$A'$' also equals $\Cl(K_n^+)/\pnpo\Cl(K_n^+)$ in this situation.
So Rubin's Theorem.~1.3 implies that $\iota_n(\{\bar{\eps}_n,\bar{h}\}_n)$ annihilates the latter (giving our
Theorem~\ref{thm: Stick in plus part mod pnpo}) provided $\eps_n$ lies in Rubin's `$\cal C$'
\ie\ it is a `special' unit.
If $f_n$ is not a power of $p$ then it is certainly a unit and the proof that it is special is similar to that of Rubin's Theorem~2.1. (Take
$u:=N_{\bbQ(\xi_{qf_n})/K_n^+\bbQ(\xi_q)^+}(1-\xi_{q}\xi_{f_n})$ for each prime  $q\neq 2$ splitting in $K_n^+$.)
If $f_n$ is a power of $p$ then $\eps_n$ is only a `special number' but see Rubin's Remark~2, p.~513.\vertsp\\
\noindent Let us define a subset $\bar{\fD}_n$ of $\cR_n[G_n]^+$ for $n\geq 0$ by
\beq\label{eq: def of Dnbar}
\bar{\fD}_n:=\left\{\iota_n(\{\bar{\eta}_n,\bar{h}\}_n):\bar{h}\in\bar{\fX}_n\right\}
=\left\{\iota_n(\{\bar{\eta}_n,\bar{h}\}_n):\bar{h}\in\bar{\fX}_n^-\right\}
\eeq
which is clearly an ideal since $\bar{h}\mapsto\iota_n(\{\bar{\eta}_n,\bar{h}\}_n)$ is
$\cR_n [G_n]$-semilinear w.r.t.\ $\iota_n$.
Since every $\bar{h}\in\bar{\fX}_n$ is a Frobenius element in $P_n/K_n$.
we can use Proposition~\ref{prop: pairing and power residue} to reformulate
Theorem~\ref{thm: Stick in plus part mod pnpo} as the following remarkable strengthening of Claim~\ref{claim: for thm 1}.
\begin{cor}\label{cor: to Stick in plus part mod pnpo} If $n\geq 0$ then
\beq\label{eq: explicit form for Lnbar}
\bar{\fD}_n=\left\{\half\sum_{g\in G_n}\left\{\frac{g\inv(\eps_n)}{\fq}\right\}_n g:\mbox{$\fq$ a prime of $K_n$ not dividing $p$}\right\}.
\eeq
Moreover, $\bar{\fD}_n$ is an ideal of $\cR_n[G_n]^+$ annihilating $\Cl(K_n^+)/\pnpo\Cl(K_n^+)$. \ePf
\end{cor}
\noindent We now pass to limits as $n\rightarrow\infty$, as
explained in Section~\ref{sec: gencon}. If $m\geq n \geq 0$ then Lemma~\ref{lemma: norm relations} implies $N^m_n \eps_{K_m}=\eps_{K_n}$
so $N^m_n \eps_m=\eps_n$.
Thus $\underline{\eta}:=(\eta_k)_{k\geq
0}$ lies in $\cV_\infty^+$ and
we may define a (continuous) $\Lambda_G$-linear map $\jinf$ by
\[
\begin{array}{rcrcl}
 \jinf &:& \Xinf\ &\longrightarrow& \Lambda_G \\
    & &    &                    \\
    & & h &\longmapsto    & \{\underline{\eta},h\}_\infty
\end{array}
\]
Since $c\underline{\eta}=\underline{\eta}$, equation~(\ref{eq: Lambda-action for equiv pairing at infty}) shows as before that $\jinf$ takes
values in $\Lambda_G^-$ and factors through the projection of $\Xinf$ on $\fX_\infty^-$. We write $\fJ_\infty$ for the (closed) ideal
$\im(\jinf)$ of $\Lambda_G^-$.

For each $n\geq 0$, we may also consider the composite map
$\phi_n^\infty\circ\fj_\infty:\fX_\infty\rightarrow\bbZ_p[G_n]$
where $\phi_n^\infty$ is the natural map
$\Lambda_G\rightarrow\bbZ_p[G_n]$. Clearly, this factors through the module
$\fX_{\infty,\Gamma_n}$ of $\Gamma_n$-covariants. Now $\Gamma_n$ is
pro-cyclic, generated by $\gamma_n$, say, so that
$\fX_{\infty,\Gamma_n}=\Xinf/(1-\gamma_n).\Xinf$ and a well-known
argument shows that $(1-\gamma_n).\Xinf$ is the (closure
of the) commutator subgroup of $\Gal(M_\infty/K_n)$, namely
$\Gal(M_\infty/M_n)$. It follows that the natural map
$\fX_{\infty,\Gamma_n}\rightarrow \fX_n$ is {\em injective} with
image $\fX_n^0:=\Gal(M_n/K_\infty)\subset\fX_n$. Therefore,
$\phi_n^\infty\circ\fj_\infty$ factors through a unique
$\bbZ_p[G_n]$-linear map
\[
\fj_n:\fX_n^0\longrightarrow \bbZ_p[G_n]
\]
Unravelling the above definitions, that of $\{\
,\ \}_\infty$ and the identification of $\Lambda_G$ with
${\displaystyle\lim_{\leftarrow}}\ \cR_n[G_n]$ in Section~\ref{sec: gencon},
we obtain the following, more explicit description of
$\fj_n(h)$ for any $h\in\fX_n^0$:
\[
\fj_n(h)=\lim_{m\rightarrow\infty}\phi^m_n(\{\bar{\eta}_m,\bar{h}_m\}_m)
\]
where, for each $m\geq n$,  $h_m$ is any lift of $h$ to
$\fX_m^0$ (the choice does not matter) and
$\phi^m_n:\cR_m[G_m]\rightarrow\cR_m[G_n]$ is as in
Section~\ref{sec: gencon}. Clearly, $\fj_n$
factors through the projection of $\fX_n^0$ on
$(\fX_n^0)^-=\fX_n^-$ and $\im(\fj_n)=\phi^\infty_n(\fJ_\infty)$ is an ideal of $\bbZ_p[G_n]^-$
which we shall denote $\fJ_n$. We shall examine $\fj_\infty$,  $\fJ_\infty$, $\fj_n$ and $\fJ_n$ more closely in Section~\ref{sec: inertia subgps}.

Now let us write $\fX_\infty^\dag$ for the module $\fX_\infty$ with the $\Lambda_G$-action twisted by $\iota_\infty$.
The composite map $\fd_\infty:=\iota_\infty\circ\jinf:\fX_\infty^\dag\rightarrow\Lambda_G$
(taking $h$ to $\iota_\infty(\{\underline{\eta},h\}_\infty)$)
is then continuous and $\Lambda_G$-linear and factors through the projection on $(\fX_\infty^-)^\dag=(\fX_\infty^\dag)^+$.
We set
$\fD_\infty:=\im(\fd_\infty)=\iota_\infty(\fJ_\infty)$, which is a (closed) ideal of $\Lambda_G^+$ and for each $n\geq 0$
we write $\fD_n$ for
$\im(\phi_n^\infty\circ\fd_\infty)$, \ie
\beq\label{eq: def. of Ln}
\fD_n=\phi^\infty_n(\fD_\infty)=\{\phi^\infty_n(\fd_\infty(h)):h\in\fX_\infty^-\}\subset\bbZ_p[G_n]^+
\eeq
Clearly, $\fD_n$ is an ideal of $\bbZ_p[G_n]^+$  and the latter will henceforth be identified with $\bbZ_p[G_n^+]$.
Since the map $h\mapsto h|_{P_n}$ sends $\fX_\infty^-$ {\em onto} $ \bar{\fX}_n^-$,
the reduction of $\fD_n$ modulo $\pnpo$ in $\cR_n[G_n]^+$ is the ideal previously denoted
$\bar{\fD}_n$.
If $\bar{\fJ}_n$ denotes the corresponding reduction of $\fJ_n$ in $\cR_n[G_n]^-$ then clearly
$\iota_n(\bar{\fJ}_n)=\bar{\fD}_n$ but there appears to be no direct relation between $\fJ_n$ and $\fD_n$ themselves.

We now give an abstract description of $\fD_n$. If $h\in\fX_\infty$ then
\begin{eqnarray}
(\mbox{coefficient of $g$ in $\phi^\infty_n(\fd_\infty(h))$})&=&
(\mbox{coefficient of $g$ in $\phi^\infty_n(\iota_\infty(\{\underline{\eta},h\}_\infty))$})\nonumber\\
&=&(\mbox{coefficient of $g$ in ${\displaystyle\lim_{m\rightarrow\infty}}\phi^m_n(\iota_m(\{\bar{\eta}_m,h|_{P_m}\}_m))$})\nonumber\\
&=&\lim_{m\rightarrow\infty}\sum_{\pi^m_n(\tg)=g }\langle \tg\inv\bar{\eta}_m,h|_{P_m}\rangle_m\nonumber\\
&=&\lim_{m\rightarrow\infty}\langle g\inv\bar{\eta}_n,h|_{P_m}\rangle_m\label{eq: characterisation of Ln}
\end{eqnarray}
using~(\ref{eq: formula for iotan applied to pairing}) and the fact $N^m_n\eta_m=\eta_n$.
Of course, {\em any} element
$\alpha$ of $\bbZ_p\otimes\cV_n$ gives an element $\bar{\alpha}\in\bar{\cV}_m$ for all $m\geq n$ and it is easy to see
for any $h\in \fX_\infty$ the limit
\[
\lfloor\alpha,h\rceil^\infty_n:=\lim_{m\rightarrow\infty}\langle \bar{\alpha},h|_{P_m}\rangle_m 
\]
is a well-defined element of $\bbZ_p$ which is $\bbZ_p$-bilinear as a function of $\alpha$ and $h$.
So~(\ref{eq: characterisation of Ln}) gives:
\beq\label{eq: use this}
\phi^\infty_n(\fd_\infty(h))=\sum_{g\in G_n}\lfloor g\inv\eta_n,h\rceil^\infty_n g
\eeq
\begin{prop}\label{prop: characterisation of Dn}
$\fD_n=\{F(\eta_n):F\in\Hom_{\bbZ_p[G_n]}((\bbZ_p\otimes\cV_n)^+,\bbZ_p[G_n])\}$ for all $n\geq 0$.
\end{prop}
\bPf\ For any $\bbZ_p[G_n]$-module $M$, there is a functorial isomorphism from
$\Hom_{\bbZ_p}(M,\bbZ_p)$ to $\Hom_{\bbZ_p[G_n]}(M,\bbZ_p[G_n])$
sending $f$ to the map $F:m\mapsto\sum_{g\in G_n}f(g\inv m)g\ \forall m\in M$. Thus
by~(\ref{eq: def. of Ln}) and~(\ref{eq: use this}), it suffices to show that any element of
$\Hom_{\bbZ_p}((\bbZ_p\otimes \cV_n)^+,\bbZ_p)$ is of form
$\alpha\mapsto\lfloor\alpha,h\rceil^\infty_n$ for some
$h$ in $\fX_\infty^-$. Using the compactness of $\fX_\infty^-$ and the definition of $\lfloor\cdot,\cdot\rceil^\infty_n$,
we are reduced to showing the surjectivity of the following composite map for all $m\geq n$:
\[
\fX_\infty^-\stackrel{a_m}{\longrightarrow}\Hom_{\cR_m}(\bar{\cV}_m^+,\cR_m)
\stackrel{b_{m,n}}{\longrightarrow}\Hom_{\cR_m}((\cV_n/\cV_n^{p^{m+1}})^+,\cR_m)
\]
where $a_m(h)$ is the homomorphism $\bar{\alpha}\mapsto\langle\bar{\alpha},h|_{P_m}\rangle_m$ and
$b_{m,n}$ is induced by the restriction of the natural map $\cV_n/\cV_n^{p^{m+1}}\rightarrow\cV_m/\cV_m^{p^{m+1}}=\bar{\cV}_m$
to plus-parts. But it is an easy exercise
to see that the latter restriction is injective and since $\cR_m$ is injective as a module over itself, it follows that $b_{m,n}$ is surjective for all
$m\geq n$. Furthermore, the surjectivity of $a_m$ for all $m\geq n$  is an immediate consequence of that of
$\fX_\infty^-\rightarrow\bar{\fX}_m^-$ and Kummer theory, taking into account the fact that $\cV_m^{p^{m+1}}=\cV_m\cap(K_m^\times)^{p^{m+1}}$.
The result follows.\ePf
\rem\label{rem: replacing by plus part} Since $\bbZ_p\otimes E_{S_p}(K_n^+)=(\bbZ_p\otimes\cV_n)^+$ is a
$\bbZ_p[G_n]$-direct summand of $\bbZ_p\otimes\cV_n$, we can of course replace
the former by the latter in the statement of the Proposition. If $f_n$ is not a power of $p$ then $\eta_n$ actually lies in
$\bbZ_p\otimes E(K_n^+)$ which is a $\bbZ_p$-direct summand of $\bbZ_p\otimes E_{S_p}(K_n^+)$. By the functoriality mentioned in the proof,
it follows that in this case we can also replace $(\bbZ_p\otimes\cV_n)^+$ by $\bbZ_p\otimes E(K_n^+)$ in the Proposition.\vertsp\\
For each prime number $r$ we write $D_r(K_n/\bbQ)$ for the common decomposition subgroup of $G_n$ at primes of $K_n$ above $r$ and
$N_{D_r(K_n/\bbQ)}$ for the norm element $\sum_dd\in\bbZ[G_n]$ where $d$ runs through $D_r(K_n/\bbQ)$. If $r|f_n$ 
and $f_n$ is not a prime power (\ie\ not a power of $p$) then Lemma~\ref{lemma: norm relations} implies
$N_{D_r(K_n/\bbQ)}\eta_n=1$ and Prop.~\ref{prop: characterisation of Dn} gives
\begin{prop}\label{prop: Ln killed by norm-elts.}
Suppose $n\geq 0$ and $f_n$ is not a power of $p$.
Then $N_{D_r(K_n/\bbQ)}\fD_n=\{0\}$   for every prime number $r$ dividing $f_n$ (\eg\ $r=p$). In particular, $\bbZ_p[G_n^+]/\fD_n$
is infinite.
\ePf
\end{prop}
Let $n_0$ be the smallest value of $n\geq 0$ such that $K_\infty/K_n$ is ramified at one (hence any) prime above $p$. Thus
$\Gamma_{n_0}$ is precisely the  inertia subgroup of
$\Gamma_0\cong\bbZ_p$ at any prime above $p$. For each $n\geq -1$ we let $L_n$ denote the maximal unramified abelian $p$-extension of $K_n$, so
$L_n\subset M_n$ and we write $L_\infty$ for $\bigcup_{n\geq -1}L_n\subset M_\infty$.
Then $X_n:=\Gal(L_n/K_n)$ and $X_\infty:=\Gal(L_\infty/K_\infty)$ are
isomorphic  via the Artin maps to $A_n:=\Cl(K_n)_p$ and ${\displaystyle \lim_\leftarrow }\,A_m$
(limit w.r.t.\ the norm maps $N^m_n$) respectively, as modules for $\bbZ_p[G_n]$ and
$\Lambda_G$.
We may identify $A_n^+$ with $\Cl(K_n^+)_p$ and consider it  as a $\bbZ_p[G_n^+]$-module. If $\cR[H]$ is any group-ring, we shall write
$I(\cR[H])$  for its augmentation ideal. We can now state our main annihilation result.
\begin{thm}\label{thm: Stick in + part over Zp v2}
Let $K$ be as above and $n\geq 0$. Then
\begin{enumerate}
\item\label{part: 1 of Stick in + part over Zp v2}
$\fD_\infty$ annihilates ${\displaystyle \lim_\leftarrow }\,A_m^+$ (or, equivalently, $X_\infty^+$).
\item\label{part: 2 of Stick in + part over Zp v2}
If $n\geq n_0$ then  $\fD_n$  annihilates $A_n^+$ (or, equivalently, $X_n^+$).
\item\label{part: 3 of Stick in + part over Zp v2}
In any case
$I(\bbZ_p[G_n^+])\fD_n$ annihilates $A_n^+$ (or, equivalently, $X_n^+$).
\end{enumerate}
\end{thm}
\bPf\ Suppose $h\in\fX_\infty^\dag$ and $(\fc_m)_{m}\in {\displaystyle\lim_\leftarrow }\,A_m^+$ and set $\fd_\infty(h)(\fc_m)_{m}=(\fb_m)_{m}$.
By definition
$\fb_m=\phi^\infty_m(\iota_\infty(\{\underline{\eta}, h\}_\infty))\fc_m$ where
$\phi^\infty_m(\iota_\infty(\{\underline{\eta}, h\}_\infty))\in \bbZ_p[G_m^+]$ is congruent modulo $p^{m+1}$ to
$\iota_m(\{\bar{\eta}_m, h|_{P_m}\}_m)$. So Theorem~\ref{thm: Stick in plus part mod pnpo} implies $\fb_m\in p^{m+1}A_m^+$ for all $m$.
Thus, for any $n\geq 0$, $\fb_n=N^m_n\fb_m\in p^{m+1}A_n^+$ for any $m\geq n$ so $\fb_n=0\ \forall n$. This proves
part~\ref{part: 1 of Stick in + part over Zp v2}. If $n\geq n_0$ then $K_\infty/K_n$ is totally ramified above $p$ so
the restriction $X_\infty\rightarrow X_n$ is surjective
and~\ref{part: 2 of Stick in + part over Zp v2}
follows from~\ref{part: 1 of Stick in + part over Zp v2}.
Part~\ref{part: 3 of Stick in + part over Zp v2} follows similarly (for $n\leq n_0$) using the fact that the cokernel of
$X_\infty\rightarrow X_n$ is $\Gal((K_\infty\cap L_n)/K_n)=\Gal(K_{n_0}/K_n$), on which $G_n$ clearly acts trivially.\ePf
\rem\ \label{rem: m_0, i_0, F etc.}To clarify the picture, define integers $i_0, f'$ and $m_0$ by $\mu_{p^\infty}(K_0)=\mu_{p^{i_0+1}}$ and $f_0=f'p^{m_0+1}$ with $p\ndiv f'$.
It is easy to see that $m_0\geq i_0\geq 0$ and that $K_0=K_1=\ldots=K_{i_0}$ while $[K_n:K_0]=p^{n-i_0}$ for $n>i_0$. One shows that
$\bbQ(\mu_{f_0})/K_0$ is
unramified above $p$, hence also $K_{m_0}/K_0$. On the other hand, if $F$ is the
inertia subfield of $K_{m_0}$ at $p$, one can show that $K_{m_0}=F_{m_0}$, so $K_\infty/K_{m_0}=F_\infty/F_{m_0}$ is totally ramified above $p$.
Hence $K_{n_0}=K_{m_0}$. Thus $n_0$ equals $m_0$ or $0$  according as $m_0>i_0$  or $m_0=i_0$, and $[K_{n_0}:K_0]=p^{m_0-i_0}$
in both cases.\vertsp\\
\rem\label{rem: Gamma_n covts}
The module $(X_\infty^+)_{\Gamma_n}=X_\infty^+/(1-\gamma_n)X_\infty^+$
is finite since it is a quotient of $(\fX_\infty^+)_{\Gamma_n}\cong\fX_n^{0,+}$ which is finite by Leopoldt's Conjecture for
$K_n^+$ (which holds \eg\ by~\cite[Thm.~5.25]{Wash}.) Suppose for simplicity that $n\geq n_0$. Then the map
$X_\infty^+\rightarrow X_n^+$
factors through a surjection $y_n^+:(X_\infty^+)_{\Gamma_n}\rightarrow X_n^+$.
Part~\ref{part: 1 of Stick in + part over Zp v2} of the Theorem clearly implies that
$\fD_n$ annihilates $(X_\infty^+)_{\Gamma_n}$  as $\bbZ_p[G_n^+]$-module, which is {\em a  priori} a stronger statement
than part~\ref{part: 2 of Stick in + part over Zp v2} whenever
$\ker(y_n^+)=:Y_n^+/(1-\gamma_n)X_\infty^+$ is non-trivial. However, one can show that $\ker(y_n^+)=\{0\}$ if $|S_p(K_n^+)|=1$ and that
$Y_{n+1}^+=(1+\gamma_n+\gamma_n^2+\ldots+\gamma_n^{p-1})Y_n^+$ in general (see \eg~\cite[Lemma 13.15]{Wash}).
The latter implies that
$\ker(y_{n+1}^+)$ is a quotient of $\ker(y_n^+)$ and hence that
$\ker(y_i^+)$ is finite and decreasing in size as $i\rightarrow\infty$. In particular, it must stabilise.
\vertsp\\
\rem\label{TESTAGAIN}
If $m\geq n\geq 0$  then following generalisation of~(\ref{eq: explicit form for Lnbar}) can be deduced
from~(\ref{eq: characterisation of Ln}) (for example).
\beq\label{eq: explicit Ln mod pmpo}
\beta_{\infty,m;n}(\fD_n)=\phi^m_n(\bar{\fD}_m)=
\left\{\half\sum_{g\in G_n}\left\{\frac{g\inv(\eps_n))}{\fq}\right\}_m g:
\mbox{$\fq$ a prime of $K_m$ not dividing $p$}
\right\}
\eeq
If also $n\geq n_0$ and $p^{m+1}$ kills $A_n^+$ then the R.H.S.\ above provides an explicit annihilator of $A_n^+$ in
$\cR_m[G_n^+]$. (If $n<n_0$ we may have to multiply by $I(\cR_m[G_n^+])$.) One would like to `let $m$ tend to infinity' and obtain
a similar expression for $\fD_n$ itself. Unfortunately, this cannot be done,
essentially because the limited splitting of finite primes in the extension $K_\infty/K_n$ means that the sequence
$\{\sigma_{\fq_{m}, P_{m}/K_{m}}\}_{m\geq n}$ for a `coherent' sequence of such primes $\fq_{m}$ (of $K_m$)  can
never cohere to give an element of $\fX_\infty$ as $m\rightarrow\infty$.
\section{Variation of $K$, the Kernels of $\fj_\infty$ and $\fd_\infty$, and Greenberg's Conjecture}\label{sec: var K and kernels}
First we compare the maps and objects defined above for an abelian field $K$ with the corresponding ones
defined identically for a subfield $F\subset K$.
To distinguish them we shall sometimes need to include the field in the notation, using a subscript for maps and parentheses for objects.
(If omitted, the field is $K$.)
For each $n\geq 0$, and also `for $n=\infty$', the fields $K_n$ and
$M_{n}(K)$ contain $F_n$ and $M_n(F)$ respectively, so we get continuous restriction maps of
Galois groups $\pi_{K_n/F_n}:G_n(K)\rightarrow G_n(F)$ and $\fr_{n,K/F}:\fX_n(K)\rightarrow\fX_n(F)$.
Extending the former by $\bbZ_p$-linearity, we get ring homomorphisms
$\pi_{K_n/F_n}:\bbZ_p[G_n(K)]\rightarrow \bbZ_p[G_n(F)]$ for every integer $n\geq 0$ whose limit is a continuous homomorphism
$\Lambda_G(K)\rightarrow\Lambda_G(F)$ extending $\pi_{K_\infty/F_\infty}$ and denoted by the same symbol.
The prime factors of $f_{F_0}$ (which include $p$) coincide with those
of $f_{F_n}$ for all $n\geq 0$. For any other prime $r$, the Frobenius elements $\sigma_{r,F_n/\bbQ}$ for $n\geq 0$
cohere to give an element
$\sigma_{r,F_\infty/\bbQ}$ of $G_\infty(F)$ satisfying $\chi_{cyc,F}(\sigma_{r,F_\infty/\bbQ})=r\in \bbZ_p^\times$.
Define elements $x_{\infty,K/F}$ and $y_{\infty,K/F}$ of $\Lambda_G(F)$ by setting
\[
x_{\infty,K/F}:=\prod_{r|f_{K_0}\atop r\ndiv f_{F_0}}(1-\sigma_{r,F_\infty/\bbQ}\inv)\ \ \ \ \mbox{and} \ \ \ \ \
y_{\infty,K/F}:=\iota_{\infty,F}(x_{K/F})=\prod_{r|f_{K_0}\atop r\ndiv f_{F_0}}(1-r\inv\sigma_{r,F_\infty/\bbQ})
\]
For any integer $n\geq 0$, their images in $\bbZ_p[G_n(F)]$ under $\phi^\infty_{n,F}$ are denoted $x_{n,K/F}$ and $y_{n,K/F}$ respectively
and are given by the same products with $\sigma_{r,F_n/\bbQ}$
replacing $\sigma_{r,F_\infty/\bbQ}$.
\begin{prop} Let $K$, $F$ and notations be as above. Then
\beq\label{eq: relations between j and d for K and F}
\pi_{K_\infty/F_\infty}\circ \fd_{\infty,K}=x_{\infty,K/F}\,\fd_{\infty,F}\circ \fr_{\infty,K/F}
\ \ \ \mbox{and}\ \ \ \pi_{K_\infty/F_\infty}\circ \fj_{\infty,K}=y_{\infty,K/F}\,\fj_{\infty,F}\circ \fr_{\infty,K/F}
\eeq
Moreover, $\fr_{n,K/F}$ maps $\fX_n(K)^-$ {\em onto} $\fX_n(F)^-$ for any $n\geq 0$ and for $n=\infty$.
\end{prop}
\bPf\ For every $n\geq 0$ the norm $N_{K_n/F_n}$ induces a map $\bar{N}_{K_n/F_n}:\bar{\cV}_n(K)\rightarrow \bar{\cV}_n(F)$ and it follows from
Lemma~\ref{lemma: norm relations} that $\bar{N}_{K_n/F_n} \bar{\eta}_{n,K}=x_{n,K/F}\bar{\eta}_{n,F}$.
If $h\in M_\infty(K)$ then equation~(\ref{eq: formula for iotan applied to pairing}) gives
\begin{eqnarray*}
\pi_{K_n/F_n}(\iota_{n,K}\{\bar{\eta}_{n,K}, h|_{P_n(K)}\}_{n,K})&=&\sum_{g'\in G_n(F)}\sum_{g\in G_n(K)\atop \pi_{K_n/F_n}(g)=g'}
\langle g\inv\bar{\eta}_{n,K},h|_{P_n(K)}\rangle_{n,K}g'\\
&=&\sum_{g'\in G_n(F)}
\langle {g'}\inv\bar{N}_{K_n/F_n}\bar{\eta}_{n,K},h|_{P_n(F)}\rangle_{n,F}g'\\
&=&\iota_{n,F}\{x_{n,K/F}\bar{\eta}_{n,F},h|_{P_n(F)}\}_{n,F}
\end{eqnarray*}
Taking the inverse limit over $n$, using the definition of $\fd_{\infty,F}$  and equation~(\ref{eq: Lambda-action for equiv pairing at infty}) then gives
\[
\pi_{K_\infty/F_\infty}(\fd_{\infty, K}(h))=
\iota_{\infty,F}\{x_{\infty,K/F}\,\underline{\eta}_F, \fr_{\infty,K/F}(h)\}_{\infty,F}=
x_{\infty,K/F}\,\fd_{F,\infty}(\fr_{\infty,K/F}(h))
\]
whence the first equation in~(\ref{eq: relations between j and d for K and F}). The second follows on applying
$\iota_{\infty,F}$ since $\iota_{\infty,F}\circ \pi_{K_\infty/F_\infty}=\pi_{K_\infty/F_\infty}\circ \iota_{\infty,K}$.
For the final statement, use the fact that $c$ acts trivially on $\coker(\fr_{n,K/F})$ since the latter is isomorphic to
$\Gal((M_n(F)\cap K_n)/F_n)$ and $K_n$ is abelian over $\bbQ$.\ePf
\noindent It is a simple exercise to deduce the
\begin{cor}\label{cor: going down for D, j and J} For each integer $n\geq 0$:
\begin{enumerate}
\item \label{part: 1 of going down for D, j and J} $\pi_{K_n/F_n}(\fD_n(K))$ equals
$x_{n,K/F}\fD_n(F) \subset \bbZ_p[G_n(F)]^+$.
\item $\pi_{K_n/F_n}\circ \fj_{n,K}$ equals $y_{n,K/F}\,\fj_{n,F}\circ \fr_{n,F/K}$
as a map $\fX_n(K)^-\rightarrow \bbZ_p[G_n(F)]^-$ and $\pi_{K_n/F_n}(\fJ_n(K))$ equals
$y_{n,K/F}\fJ_n(F)\subset \bbZ_p[G_n(F)]^-$.  \ePf
\end{enumerate}
\end{cor}
Note that part~\ref{part: 1 of going down for D, j and J} also follows from Prop.~\ref{prop: characterisation of Dn} and implies
Prop.~\ref{prop: Ln killed by norm-elts.} for $r\neq p$. (Take $F$ to be the splitting field of $r$ in $K$ \etc)\vertsp\\
\rem\label{rem: fattening D} In contrast to Prop.~\ref{prop: Ln
killed by norm-elts.}, we shall see later  (Cor.~\ref{cor: fJn
contains fSn}) that the index of $\fJ_n(K)$ in $\bbZ_p[G_n(K)]^-$ is
often  (perhaps always) finite for all $n\geq 0$. Nevertheless,
Corollary~\ref{cor: going down for D, j and J}
suggests the
idea of `enlarging' both $\fD_n(K)$ and $\fJ_n(K)$ by a
method similar to that often used for the Stickelberger ideal
(see \eg~\cite[\S2]{Greither}): one should add to each of
them ($\bbZ_p$-multiples of) the images under $\cores_{F_n}^{K_n}$
of the corresponding ideals for all subfields $F$ of $K$. Here
$\cores_{F_n}^{K_n}$ denotes the additive homomorphism from
$\bbZ_p[G_n(F)]$ to $\bbZ_p[G_n(K)]$ which sends $g\in G_n(F)$ to the
sum of its pre-images in $G_n(K)$ under restriction. The same can
be done at the infinite level. (Indeed, if $K/F$ is not linearly disjoint from $F_\infty/F$,
one should do this first, then get the `enlarged' ideals at finite levels as images under
$\phi^\infty_n$.)
Annihilation statements similar to those of Theorem~\ref{thm: Stick in + part over Zp v2}
can be proven for the enlarged ideals $\tilde{\fD}$ using the original
versions for all subfields $F$ of $K$.
\vertsp\\
For any $F\subset K$ we write $E^0(F)$ for the subgroup of $E(F)$  consisting of those units whose local absolute norms are
trivial at all primes dividing our fixed prime $p$. Since $F$ is abelian, this is simply the kernel of $N_{D_p(F/\bbQ)}$ acting on
$E(F)$. Now let $N_\infty=N_\infty(K)$ (\resp\ $N_\infty^0=N_\infty^0(K)$) denote the infinite abelian extension of $K_\infty$ obtained by adjoining to it all $p$-power
roots of all elements of $E(K_n^+)$ (\resp\ of $E^0(K_n^+)$) for all $n\geq 0$.
Both $N^0_\infty$ and $N_\infty$ are Galois over $\bbQ$ and it is easy to see that
$N^0_\infty\subset N_\infty\subset M_\infty^-$. (Here, $M_\infty^-$ is defined by $\Gal(M_\infty/M_\infty^-)=\fX_\infty^+$ so that
$\fX_\infty^-$ maps isomorphically onto $\Gal(M_\infty^-/K_\infty)$.)
Since $K_n$ is CM it is well known that
$|E(K_n):\mu(K_n)E(K_n^+)|=1$ or $2$. It follows that we could have used $E(K_n)$  in place of
$E(K_n^+)$ (\resp\ $E^0(K_n)$ in place of $E^0(K_n^+)$) in the definition of $N_\infty$ (\resp\ of $N^0_\infty$).
\begin{thm}\label{thm: ker of jinf and dinf}
$\ker(\fd_{\infty})=\ker(\fj_{\infty})=\Gal(M_\infty/N^0_\infty)$ as subgroups of $\fX_\infty$.
\end{thm}
Before giving the proof, we deduce a first link with Greenberg's Conjecture for the extension $K_\infty^+/K^+$, \ie\ the statement
that  $|A_n^+|$ is bounded as $n\rightarrow\infty$ or, equivalently, that $X_\infty^+$ is finite.
Proposition~2 of~\cite{Greenberg} shows that
this is also equivalent to the triviality of
$A_\infty^+$ where $A_\infty$  denotes
{\em direct} limit of the $A_n$'s as $n\rightarrow\infty$ w.r.t.\ the maps coming from extension of ideals.
Now, Kummer Theory gives a non-degenerate, Galois-equivariant pairing
\beq\label{eq: gal-equiv pairing}
\Gal(M_\infty/N_\infty)\times A_\infty\rightarrow\mu_{p^\infty}
\eeq
(This follows from~\cite{Wash} pp. 294-295.
See also \cite[Thm.\ 14]{Iwasawa on Z_l extensions...} using roots of $p$-units and $p$-class groups instead.)
Hence Greenberg's Conjecture for $K_\infty^+/K^+$ is also equivalent to $\Gal(M_\infty/N_\infty)^-=\{0\}$ \ie\ $M_\infty^-=N_\infty$.
\begin{cor}\label{cor: what happens when Sp has one element I}
Suppose $|S_p(K_{n_0}^+)|=1$. Then $N_\infty^0=N_\infty$ and
\beq\label{eq: kernels with only one prime above p}
\ker(\fd_{\infty})\cap(\fX_\infty^\dag)^+=
(\ker(\fj_{\infty})\cap\fX_\infty^-)^\dag=
(\Gal(M_\infty/N_\infty)^-)^\dag\cong
 \Hom_{\bbZ_p}\left(A^+_\infty,\bbQ_p/\bbZ_p\right)
\eeq
as $\Lambda_G^+$-modules.
In particular, Greenberg's Conjecture holds for $K_\infty^+/K^+$ if and only if $\fd_\infty$ is injective on $(\fX_\infty^\dag)^+$
(or, equivalently, $\fj_\infty$ on $\fX_\infty^-$).
\end{cor}
\bPf\ By total ramification,
$|S_p(K_{n_0}^+)|=1$ is equivalent to $|S_p(K_n^+)|=1$  for all $n\geq 0$. This implies $D_p(K_n^+/\bbQ)=G_n^+$ so that
$E(K_n^+)^2\subset E^0(K_n^+)$ for all $n\geq 0$, and hence $N_\infty^0=N_\infty$.
The second equality now follows from the Theorem, as does the first (since $\fX_\infty^-=(\fX_\infty^\dag)^+$ as groups).
The isomorphism follows from the above pairing and the rest is a direct consequence.\ePf
\noindent In Section~\ref{sec: inertia subgps}, equation~(\ref{eq: double iso}) will show
that $\ker(\fj_\infty)\cap \fX_\infty^-= \Gal(M_\infty^-/N_\infty^0)$ contains
a specific submodule  which is non-trivial (and infinite) whenever $|S_p(K_{n_0}^+)|>1$, regardless of Greenberg's Conjecture.
Nevertheless, Theorem~{\ref{thm: Lambda_Gamma Structure}}\ref{part: 2 of Lambda_Gamma Structure} will show
that $\ker(\fj_\infty)$ is still as small as possible for a map $\fX_\infty\rightarrow\Lambda_G$.\vertsp\\
{\sc Proof of Theorem~\ref{thm: ker of jinf and dinf}}. The first equality follows from the injectivity
of the involution $\iota_\infty$. For any abelian field $F$ we temporarily denote by $B_\infty(F)$ the fixed field of
$\ker(\fj_{\infty,F})$ acting on $M_\infty(F)$.
The second equality thus amounts to $B_\infty(K)=N_\infty^0(K)$ which we now prove.
\begin{lemma}\label{lemma: first description of Binfty}
Let $F$ be any abelian field, $n\geq 0$ and write $E_{S_p}(F_n^+)$
as a left $\bbZ[G_n(F)]$-module (under multiplication). The field
$B_\infty(F)$ is obtained by adjoining to $F_\infty$ the following subset of $M_\infty(F)^-$
\[
\cS_\infty(F):=\{\alpha^{1/p^m}:\alpha\in I(\bbZ[G_n(F)])\eps_n(F)\ \mbox{and}\ m,n\geq 0 \}
\]
Furthermore, $I(\bbZ[G_n(F)])\eps_n(F)$ is contained in $E^0(F_n^+)$.
\end{lemma}
\bPf\ If $h\in \fX_{\infty, F}$, then $\fj_{\infty,F}(h)=0$ iff $\{\bar{\eta}_{n}(F), h|_{P_n(F)}\}_{n,F}=0$ for all $n\geq 0$. This is equivalent to
$h$ fixing all conjugates of $\eps_n(F)^{1/\pnpo}$ over $\bbQ$ for all $n\geq 0$
and hence to $h$ fixing all conjugates of $\eps_n(F)^{1/p^m}$ for all $m,n\geq 0$, since
$\eps_n(F)=N_{\Gal(F_m/F_n)}\eps_m(F)$ for $m\geq n$. This shows that $B_\infty(F)$ is obtained by adjoining the larger set with
$I(\bbZ[G_n(F)])$ replaced by $\bbZ[G_n(F)]$ in the definition of $\cS_\infty(F)$. The fact that this gives the same field follows from the relation
$\eps_n(F)\equiv (N_{\Gal(F_{m+n}/F_n)}-p^m)\eps_{m+n}(F)\pmod{(K_{m+n}^\times)^{p^m}}$ (whenever $m\geq 0$)
since $N_{\Gal(F_{m+n}/F_n)}-p^m\in I(\bbZ[G_{m+n}(F)])$.
This proves the first statement. For the second, if  $f_{F_n}$ is a power
of $p$ then $|S_p(F_n^+)|=1$, so it follows from $\eps_n(F)\in E_{S_p}(F_n^+)$. Otherwise,
Lemma~\ref{lemma: norm relations} gives $\eps_{F_n}\in E(F_n)$  and easily implies
$N_{D_p(F_n/\bbQ)}\eps_{F_n}=1$. Hence $\eps_n(F)\in E^0(F_n^+)$. \ePf
\noindent  The norm relations mean that the index $|E^0(K_n^+):I(\bbZ[G_n(K)])\eps_n(K)|$ is usually infinite for every $n\geq 0$. So the equality
$B_\infty(K)=N_\infty^0(K)$ does not follow immediately from the above Lemma. We also need the less obvious
\begin{lemma}\label{lemma: BK contains BF for all F}
Suppose $F$ is a subfield of $K_n$ for some $n\geq 0$. Then $B_\infty(K)$ contains $B_\infty(F)$.
\end{lemma}
\bPf\ First, the map $\fj_{\infty,K_n}:\fX_\infty(K_n)\rightarrow\Lambda_G(K_n)$ is formally identical to
$\fj_{\infty,K}:\fX_\infty(K)\rightarrow\Lambda_G(K)$.
Thus $B_\infty(K)=B_\infty(K_n)$ and, replacing $K$ by $K_n$, we may reduce to the case $F\subset K$. We need
to prove that $h\in \ker(\fj_{\infty,K})$ implies that $\fr_{\infty,K/F}(h)$ fixes $B_\infty(F)$, \ie\ that
$\fj_{\infty,F}\circ \fr_{\infty,K/F}(h)=0$. But if $h\in \ker(\fj_{\infty,K})$ then (\ref{eq: relations between j and d for K and F})
implies $y_{\infty,K/F}\,\fj_{\infty,F}\circ \fr_{\infty,K/F}(h)=0$
so it suffices to show that $y_{\infty,K/F}$ is a non-zero-divisor of $\Lambda_G(F)$.
This follows from the fact that its image $y_{n,K/F}$ is a non-zero-divisor of $\bbZ_p[G_n(F)]$ for all $n\geq 0$. ({\em E.g.}\ because
it divides $\prod_r(1-r^{-a_{r,n}})\in \bbZ_p\setminus\{0\}$ where $a_{r,n}\geq 1$ is the order of $\sigma_{r,F_n/\bbQ}$ in $G_n(F)$.) \ePf
\noindent
By Lemma~\ref{lemma: first description of Binfty},
the following defines a $\bbZ[G_n(K)]$-submodule of $E^0(K_n^+)$ for each $n\geq 0$
\[
C^0_n(K):=\sum_{F\subset K_n}I(\bbZ[G_n(F)])\eps_n(F)
\]
(where $F$ ranges over the (finitely many) subfields of $K_n$).
Now consider the subfield $K_\infty(\tilde{\cS}_\infty(K))$ of $N_\infty^0(K)$ obtained by adjoining to $K_\infty$ the set
\[
\tilde{\cS}_\infty(K):=\{\beta^{1/p^m}:\beta\in C^0_n(K)\ \mbox{and}\ m,n\geq 0 \}
\]
Lemma~\ref{lemma: first description of Binfty} for $F=K$ and the obvious containment $\cS_\infty(K)\subset \tilde{\cS}_\infty(K)$ imply
$B_\infty(K)\subset K_\infty(\tilde{\cS}_\infty(K))$. The reverse inclusion follows from the fact that
every element of $\tilde{\cS}_\infty(K)$ is a product of elements
of the sets $\cS_\infty(F)\subset B_\infty(F)$ for $F\subset K_n$, and hence lies in $B_\infty(K)$,
by Lemma~\ref{lemma: BK contains BF for all F}. Thus $K_\infty(\tilde{\cS}_\infty(K))=B_\infty(K)$ and to prove Theorem~\ref{thm: ker of jinf and dinf}
it only remains to show that the inclusion $K_\infty(\tilde{\cS}_\infty(K))\subset N_\infty^0(K)$ is an equality. But in view of the definitions
of $\tilde{\cS}_\infty(K)$ and $N_\infty^0(K)$, this is an easy deduction from
\begin{lemma} The index $|E^0(K_n^+):C^0_n(K)|$ is finite for all $n\geq 0$.
\end{lemma}
\bPf\ This is a variant of the proof that the (full) group of cyclotomic units of an abelian field is of finite index in its unit group,
so we leave out some
details. It suffices to show that the natural inclusion of $\bbC\otimes C^0_n(K)$ in
$\bbC\otimes E^0(K_n^+)$
is an equality. We consider these as nested $\bbC[G_n]$-submodules of $\bbC\otimes E(K_n^+)$ and show that
$r_\chi:=\dim_{\bbC}(e_\chi(\bbC\otimes C^0_n(K)))$ is at least $r'_\chi:=\dim_{\bbC}(e_\chi(\bbC\otimes E^0(K_n^+)))$ for every (irreducible)
complex character $\chi$ of $G_n$ (whose idempotent in $\bbC[G_n]$ is $e_\chi$).
Dirichlet's Theorem shows that the map
\displaymapdef{\lambda}{E(K_n^+)}{\bbC[G_n]}{\eps}{\sum_{g\in G_n}\log|g(\eps)|g\inv}
extends to a $\bbC[G_n]$-isomorphism $\lambda_\bbC:\bbC\otimes E(K_n^+)\rightarrow I(\bbC[G_n])^+$. Since
$\bbC\otimes E^0(K_n^+)$ is the kernel of $N_{D_p(K_n/\bbQ)}$ acting on
$\bbC\otimes E(K_n^+)$, it follows that $r'_\chi=0$ unless $\chi(c)=1$ and
$\chi(D_p(K_n/\bbQ))\neq\{1\}$ in which case $r'_\chi=1$. So it suffices to show that $e_\chi(\bbC\otimes C^0_n(K))\neq\{0\}$ in this latter case. Fix
such a $\chi$ and
let $F$ be the (real) subfield of $K_n$ fixed by $\ker(\chi)$ so $p$ does not split completely in $F$.
We can also regard $\chi$ as an even, non-trivial Dirichlet character modulo its conductor
$f=f_{F}$ such that $\chi(\bar{p})\neq 1$ if $p\ndiv f$.
Choose also $g_0\in G_n$ such that $\chi(g_0)\neq 1$ so that $z_\chi:=e_\chi(1\otimes (1-g_0)\eps_n(F))$ lies in $e_\chi(\bbC\otimes C^0_n(K))$. If $p|f$ then
Lemma~\ref{lemma: norm relations} shows that $N_{F_n^+/F}\eps_n(F)=\eps_{F}$ and a calculation gives
\[
\lambda_\bbC(z_\chi)=
[K_n:F_n^+](1-\chi(g_0))\sum_{a=1\atop (a,f)=1}^{f}\log|1-\xi_{f}^a|\chi\inv(\bar{a})
=[K_n:F_n^+](\chi(g_0)-1)\tau(\chi\inv)L(1,\chi)
\]
where $L(s,\chi)$ is the complex (primitive) $L$-function and $\tau(\chi\inv)$ is the Gauss Sum attached to $\chi\inv$.
(See, for example, Theorem 4.9 and the preceding pages in~\cite{Wash}.)
If $p\ndiv f$ then  $N_{F_n^+/F}\eps_n(F)=(1-\sigma\inv_{p,F/\bbQ})\eps_{F}$, giving an extra factor of $(1-\chi\inv(\bar{p}))$ in the
second two members above. In either case the third member is a product of nonzero terms, so $z_\chi\neq 0$.\ePf
\noindent This completes the proof of Theorem~\ref{thm: ker of jinf and dinf}.\ePf
\section{The Case $K=\bbQ$ : the Ideals $\fD_n$ and the Map $\fd_\infty$}\label{sec: d and D for K=Q}
If $K=\bbQ$ then $K_n=\bbQ(\mu_\pnpo)$, $f_n=\pnpo$, $n_0=i_0=0$ and $K_n$ (\resp\ $K_n^+$) has a
unique prime ideal dividing $p$, generated by
$1-\zeta_n$ (\resp\ by $\eps_n$). We abbreviate $E_{S_p}(K_n^+)$ and $E(K_n^+)$ to $\tilde{E}_n$ and $E_n$ respectively.
We write $\tilde{C}_n$ for the $\bbZ[G_n^+]$-submodule of $\tEn$ generated by $\eps_n$ and  $C_n$ for $\tCn\cap E_n$. (This is  the group of
cyclotomic units of $K_n^+$ and coincides with the group $C_n^0(K)$ defined above, in this case.) Since $\tEn=\tCn E_n$, the natural map $E_n/C_n\rightarrow \tEn/\tCn$ is an isomorphism.
It follows from \eg\
Theorem~8.2 of~\cite{Wash} that  $\tEn/\tCn$ is finite, of cardinality a power of $2$ times $|\Cl(K_n^+)|$.
In particular, $\tC_n$ has the same $\bbZ$-rank as $\tE_n$, namely $[K_n^+:\bbQ]$ so that $\tCn$
is $\bbZ[G_n^+]$-{\em free} with basis $\{\eps_n\}$.

Since $p$ is odd, $\Zpten\tEn$ is $\bbZ_p$-torsionfree so may be regarded as a submodule of $\bbQ_p\otimes\tEn$. We may also regard $\bbZ_p\otimes\tCn$ as a
$\bbZ_p[G_n^+]$-free submodule with basis $\{\eta_n\}$ and spanning $\bbQ_p\otimes\tEn$ over $\bbQ_p$.
It follows that there exists a a (unique) fractional ideal $J_n$ of $\bbQ_p[G_n^+]$ (by which we mean a $\bbZ_p[G_n^+]$-submodule of $\bbZ_p$-rank equal to $|G_n^+|$)
such that the map $J_n\rightarrow\Zpten\tEn$ sending $j$ to $j\eta_n$ is an isomorphism.
For each $x\in\bbZ_p[G_n^+]$ we let $t_{n,x}\in
\Hom_{\bbZ_p}\left((\Zpten\tEn)/(\Zpten\tCn),\bbQ_p/\bbZ_p\right)$
be the homomorphism
sending the class of $j\eta_n$ to that of the coefficient of $1$ in $xj$, for all $j\in J_n$. If $H$ is any abelian group
and $M$ any $\bbZ_p[H]$-module we shall sometimes write the Pontrjagin dual
$\Hom_{\bbZ_p}\left(M,\bbQ_p/\bbZ_p\right)$  as $M^\vee$ for brevity. We emphasise that it is always endowed
with the $\bbZ_p[H]$-action for which $h.f$ is the homomorphism $f\circ h$
({\em not} $f\circ h\inv$ as in~\cite{Kraft-Schoof} \etc) for any $h\in H$ and $f\in M^\vee$.
\begin{thm}\label{prop: gp ring mod fD iso to E/C dual} Suppose $K=\bbQ$,  $n\geq 0$  and notations are as above. The map
\displaymapdef{t_n}{\bbZ_p[G_n^+]/\fD_n}
{\Hom_{\bbZ_p}\left((\Zpten\tEn)/(\Zpten\tCn),\bbQ_p/\bbZ_p\right)}
{x\bmod \fD_n}{t_{n,x}}
is a well-defined isomorphism of $\bbZ_p[G_n^+]$-modules.
\end{thm}
\bPf\ The injection $(\Zpten\tEn)/(\Zpten\tCn)\rightarrow(\bbQ_p[G_n^+]/\bbZ_p[G_n^+])$ sending the class of $j\eta_n$
to that of $j$ (for $j\in J_n$) induces a surjection from $(\bbQ_p[G_n^+]/\bbZ_p[G_n^+])^\vee$ to
$((\Zpten\tEn)/(\Zpten\tCn))^\vee$. On the other hand, it is easy to see that every element of $(\bbQ_p[G_n^+]/\bbZ_p[G_n^+])^\vee$
sends the class of $y\in\bbQ_p[G_n^+]$ to that of the coefficient of $1$ in $xy$ for
some fixed $x\in\bbZ_p[G_n^+]$. It follows that
the map $\tilde{t}_n$ from $\bbZ_p[G_n^+]$ to $((\Zpten\tEn)/(\Zpten\tCn))^\vee$ sending $x$ to $t_{n,x}$ is surjective.
It is easy to check that
$\tilde{t}_n$ is $\bbZ_p[G_n^+]$-linear so it only remains to prove that $\fD_n=\ker(\tilde{t}_n)$.
But $\ker(\tilde{t}_n)$ is precisely the set
$\{x\in\bbZ_p[G_n^+]: xj\in\bbZ_p[G_n^+]\ \forall j\in J_n\}$. It follows easily that
$\Hom_{\bbZ_p[G_n^+]}(\bbZ_p\otimes\tEn,\bbZ_p[G_n^+])=\Hom_{\bbZ_p[G_n]}((\bbZ_p\otimes\cV_n)^+,\bbZ_p[G_n])$ is precisely
the set of maps $j\eta_n\mapsto xj$ for $x\in\ker(\tilde{t}_n)$. Taking $j=1$,
Proposition~\ref{prop: characterisation of Dn} implies $\fD_n=\ker(\tilde{t}_n)$, as required.\ePf
\noindent From the above -- and elementary properties of duals \etc\ -- we deduce:
\begin{cor}\label{prop: K=Q, Fittings, duals etc.} If $K=\bbQ$ and $n\geq 0$ then
\begin{enumerate}
\item\label{part: K=Q, Fittings, duals etc. 1}
$\bbZ_p[G_n^+]/\fD_n$ is finite and $|\bbZ_p[G_n^+]/\fD_n|=|\bbZ_p\otimes(\tEn/\tCn))|=|A_n^+|$.

\item\label{part: annihilators} $\fD_n$ is precisely the $\bbZ_p[G_n^+]$-annihilator
of $(\bbZ_p\otimes(\tEn/\tCn))^\vee\cong(\bbZ_p\otimes(E_n/C_n))^\vee$, hence
also of $\bbZ_p\otimes(\tEn/\tCn)\cong\bbZ_p\otimes(E_n/C_n)$.

\item\label{part: Fitting} $\fD_n$ is precisely the (initial) $\bbZ_p[G_n^+]$-Fitting ideal
of $(\bbZ_p\otimes(\tEn/\tCn))^\vee\cong(\bbZ_p\otimes(E_n/C_n))^\vee$. \ePf
\end{enumerate}
\end{cor}
\rem\label{rem: more on Fitt and Ann} Part~\ref{part: annihilators} above combines with Theorem~\ref{thm: Stick in + part over Zp v2}~\ref{part: 2 of Stick in + part over Zp v2}
to show that
\beq\label{eq: containment of annihilators}
\Ann_{\bbZ_p[G_n^+]}(\bbZ_p\otimes(E_n/C_n))\subset\Ann_{\bbZ_p[G_n^+]}(A_n^+)
\eeq
This may be compared with statement of Thaine's theorem
in~\cite[Thm 15.2]{Wash} (as well as the general results of~\cite{Rubin-GUICG} already cited).
The former is essentially~(\ref{eq: containment of annihilators}) generalised to allow any real abelian field $F$ in
place of $K_n^+$ {\em but} also restricted to $p$ (possibly 2) not dividing $[F:\bbQ]$.
Since $G_n^+$ is cyclic, Fitting ideals of $\bbZ_p[G_n^+]$-modules and their duals coincide
(this follows from point  Prop.~1 and point 4 of~\cite[Appendix]{Mazur-Wiles}) so we have the interesting equalities
\begin{equation}\label{eq: Fitting}
\fD_n={\rm Fitt}_{\bbZ_p[G_n^+]}(\bbZ_p\otimes(E_n/C_n))={\rm Fitt}_{\bbZ_p[G_n^+]}(A_n^+)
\end{equation}
where the first follows from part~\ref{part: Fitting} above and the second from~\cite[Thm. 1]{Greither-Cornacchia}.

It is not clear to the author whether to expect generalisations of~(\ref{eq: containment of annihilators}) and/or the equalities between
each pair of the three members in~(\ref{eq: Fitting}), when $E_n$ is replaced by
$E(F)$ for arbitrary real, abelian $F$ and $C_n$ by a suitably defined group of cyclotomic units $C(F)$.
(Before even considering first equality in~(\ref{eq: Fitting})
one would have to enlarge $\fD_n$, perhaps as in Remark~\ref{rem: fattening D}.)
However, our approach certainly suggests that it might be more natural to consider
the {\em Pontrjagin  dual} of $\bbZ_p\otimes(E(F)/C(F))$. This might even be necessary in~(\ref{eq: Fitting})
when $\Gal(F/\bbQ)$ is not $p$-cyclic.
Note also that the case $p\ndiv [F:\bbQ]$ may not be indicative here. Not only is
$\bbZ_p\otimes(E(F)/C(F))$ then $\bbZ_p[\Gal(F/\bbQ)]$-isomorphic to its dual but
its Fitting ideal and annihilator coincide since $\bbZ_p\otimes E(F)$ is cyclic over
$\bbZ_p[\Gal(F/\bbQ)]$ in this case.
\vertsp\\
\noindent\ It is easy to check that the following diagram commutes for all $m\geq n\geq 0$:
\[
 \xymatrix{
& \bbZ_p[G_m^+]/\fD_m\ar[d]\ar[rrr]^-{\sim}_-{-}&&&\Hom_{\bbZ_p}\left(\Zpten(\tE_m/\tC_m),\bbQ_p/\bbZ_p\right)\ar[d]\\
& \bbZ_p[G_n^+]/\fD_n\ar[rrr]^-{\sim}_-{-}&&&\Hom_{\bbZ_p}\left(\Zpten(\tE_n/\tC_n),\bbQ_p/\bbZ_p\right)\\
}
\]
Here, the horizontal isomorphisms are (essentially) $t_m$ and $t_n$, the left-hand vertical map is the natural surjection and the
right-hand map is induced by the natural map
$(\tEn/\tCn)\rightarrow(\tE_m/\tC_m)$.
It follows that the  transition maps on the R.H.S.\ are also surjections. (This can also be seen the injectivity of
$(\tEn/\tCn)\rightarrow(\tE_m/\tC_m)$ which follows in turn from the $\bbZ[G_m^+]$-freeness of $\tC_m$.)
Passing to inverse limits, we obtain continuous $\Lambda_G^+$-isomorphisms
\beq\label{eq: iso for coker dinfty}
\Lambda_G^+/\fD_\infty\cong\lim_{\longleftarrow\atop n\geq 0}( \bbZ_p[G_n^+]/\fD_n)\cong\lim_{\longleftarrow\atop n\geq 0}((\Zpten(\tE_n/\tC_n))^\vee)
\eeq
(The first follows from compactness arguments.) Finally,
the R.H.S.\ of~(\ref{eq: iso for coker dinfty}) is easily seen to be isomorphic to
$\Hom_{\bbZ}\left(\tE(K_\infty^+)/\tC(K_\infty^+),\bbQ_p/\bbZ_p\right)$
where $\tE(K_\infty^+):=\bigcup_{n\geq 0}\tEn$ and
$\tC(K_\infty^+):=\bigcup_{n\geq 0}\tCn$.\vertsp\\
\rem\label{rem: correspondence with K-S} We point out the connections mentioned in the Introduction
between our results and those of~\cite{Kraft-Schoof}.
Let $K$ be real quadratic so that $G_n^+=G_0^+\times \Gal(K_n^+/K_0^+)$ and let $\chi$ denote
non-trivial character of $\Gal(K/\bbQ)$ inflated to $G_0^+$.
The aim of Kraft and Schoof is to study modules for such $K$ which they denote `$A_n$' and `$C_n$' and are essentially
the $\chi$-components of  our $A_n^+$ and $(\bbZ_p\otimes (E(K_n^+)/C(K_n^+)))^\vee$ respectively.
Now, using our Cor.~\ref{cor: to Stick in plus part mod pnpo}, one can obtain a precise relation
between the $\chi$-component of our $\bar{\fD}_n$ and
the denominator on the R.H.S.\ of the last equation on p.~141 of~\cite{Kraft-Schoof}, with $k=n+1$. (Use the description of
their~$f_r$ given on p.~144, not the vaguer one on p.~140.) Combining this relation
with~\cite[Prop. 2.5]{Kraft-Schoof}
gives a sort of mod-$p^k$ analogue of our Prop.~\ref{prop: characterisation of Dn}. The above-mentioned
equation itself may be compared with our Thm.~\ref{prop: gp ring mod fD iso to E/C dual} modulo $p^k$ (where, of course, $K=\bbQ$) and
the first statement of~\cite[Theorem 2.4]{Kraft-Schoof} with our equation~(\ref{eq: iso for coker dinfty}). The ideal $I$ in this statement is
essentially the $\chi$-component of our $\fD_\infty$. Note that the Galois action on Pontrjagin duals defined  in~\cite{Kraft-Schoof}
must be changed to ours to make it consistent with their own identification at the end of the proof
of~\cite[Thm.~2.4]{Kraft-Schoof}.
\vertsp\\
Vandiver's Conjecture for $p$ states that $A^+_0=\{0\}$ or, equivalently, $A^+_n=\{0\}\ \forall\, n\geq 0$.
(For the non-trivial implication, use~\cite[Thm.~10.4]{Wash}
.) Thus, Vandiver's Conjecture strengthens  Greenberg's Conjecture for  $K_\infty^+/K^+$
and each corresponds to certain properties of $\fd_\infty$:
\begin{prop} \label{prop: vandi and green} Suppose $K=\bbQ$.
\begin{enumerate}
\item \label{part: vandi and green (i)} The following are equivalent:
\begin{enumerate}
\item  Greenberg's Conjecture holds for $K_\infty^+/K^+$,
\item $\fd_\infty$ is injective on $(\fX_\infty^\dag)^+$,
\item $\coker(\fd_\infty)$ (\ie\ $\Lambda_G^+/\fD_\infty$) is finite.
\end{enumerate}
\item \label{part: vandi and green (ii)} The following are equivalent:
\begin{enumerate}
\item  Vandiver's Conjecture holds for $p$,
\item $\fd_\infty$ is an isomorphism from $(\fX_\infty^\dag)^+$ to $\Lambda_G^+$,
\item $\fd_\infty$ is surjective (\ie\ $\fD_\infty=\Lambda_G^+$).
\end{enumerate}
\end{enumerate}
\end{prop}
\bPf\ In part~\ref{part: vandi and green (i)}, the equivalence
$(a)\Leftrightarrow (b)$ is from Cor.~\ref{cor: what happens when Sp has one element I}. For
$(a)\Leftrightarrow (c)$, Greenberg's conjecture is equivalent to the boundedness of $|\bbZ_p[G_n^+]/\fD_n|$ by
Cor.~\ref{prop: K=Q, Fittings, duals etc.}~\ref{part: K=Q, Fittings, duals etc. 1}. Now use~(\ref{eq: iso for coker dinfty})
noting that the transition maps in the limits are surjective. The equivalence $(a)\Leftrightarrow (c)$ of part~\ref{part: vandi and green (ii)}
is proved in a similar way. The implication $(b)\Rightarrow (c)$ is trivial in~\ref{part: vandi and green (ii)}
and $(c)\Rightarrow (b)$ follows from the same implication
in~\ref{part: vandi and green (i)}.\ePf
\noindent \rem\ The argument $(a)\Leftrightarrow(c)$ in~\ref{part: vandi and green (i)} shows that if  Greenberg's Conjecture holds then
$|\Lambda_G^+/\fD_\infty|=|X_\infty^+|$ (since then $X_\infty^+\cong A_n^+$ for all $n>>0$).\vertsp\\
\noindent The kernel and cokernel of $\fd_\infty$ can also be related without Greenberg's Conjecture, using instead
the Main `Conjecture' of Iwasawa Theory over $\bbQ$ (a Theorem, of course!). See Thm.~\ref{prop: consequence of Main Conj}.

\section{Inertia Subgroups, the Map $\fs_\infty$ and $\Lambda_\Gamma$-torsion}\label{sec: inertia subgps}
When $\fj_\infty$ is restricted to the product of inertia subgroups in $\fX_\infty$ we shall see that it is given by
a limit of certain rather explicit $p$-adic maps $\fs_n$ as $n\rightarrow\infty$. These are specialisations of the map
$\fs_{F/k,S}$ defined in~\cite[\S 2.4]{laser} for any abelian extension
of number fields $F/k$ with $F$ of CM-type and $k$ totally real, and for any finite set $S$
of places of $k$ containing $S^0(F/k)$ (\ie\ the infinite ones
and those ramified in $F$). We start by giving the particularly simple
definition of $\fs_{F/k,S}$ in the relevant case, namely $k=\bbQ$ and $F$ any imaginary abelian field.


For each irreducible, complex character $\chi\in\widehat{\Gal(F/\bbQ)}$ we define the
$S$-truncated  $L$-function $L_{F/\bbQ, S}(s,\chi)$
to be the unique meromorphic continuation of the function defined for $\Re(s)>1$ by the Euler product
$\prod_{q\nin S}\left(1-q^{-s}\chi(\sigma_{q,F/\bbQ})\right)\inv$.
Let
\[
a^-_{F/\bbQ,S}:=\frac{i}{\pi}\sum_{\chi\in\widehat{\Gal(F/\bbQ)}\atop\ \chi\ odd}L_{F/\bbQ, S}(1,\chi)e_{\chi\inv}\in \bbC[\Gal(F/\bbQ)]^-
\]
where $e_{\chi\inv}$ denotes the idempotent $|\Gal(F/\bbQ)|\inv\sum_{g\in\Gal(F/\bbQ)}\chi\inv(g)g\inv$ of $\bbC[\Gal(F/\bbQ)]$.
We shall also write
$a^{-,\ast}_{F/\bbQ,S}$ for the image of $a^-_{F/\bbQ,S}$ under the $\bbC$-linear involution of $\bbC[\Gal(F/\bbQ)]$ sending $g\in \Gal(F/\bbQ)$ to $g\inv$.

For {\em any} integer $l\geq 1$ we
write $G(l)$ for $\Gal(\bbQ(\mu_l)/\bbQ)=
\{\sigma_{a,l}:(a,l)=1\}$ where $\sigma_{a,l}(\xi_l)=\xi_l^a$.
If $l\geq 3$ we may take $F=\bbQ(\mu_l)$ and
$S$ to be $S_l=:\{\infty\}\cup S_l(\bbQ)$. In this case we record here
(for use in Section~\ref{sec: analysis of exact sequence}) a relatively
simple {\em `equivariant functional equation'} relating
$a^{-,\ast}_{\bbQ(\mu_l)/\bbQ,S_l}$ to the Stickelberger elements  defined for any integer $r>1$ by
\beq\label{eq: def of theta}
\theta_{\bbQ(\mu_r)/\bbQ, S_r}:=\sum_{\chi\in\widehat{G(r)}}L_{\bbQ(\mu_r)/\bbQ, S_r}(0,\chi)e_{\chi\inv}=
-\sum_{a=1\atop (a,r)=1}^r\left(\frac{a}{r}-\half\right)\sigma_{a,r}\inv\in \bbQ[G(r)]^-
\eeq
(For the second equality above see \eg~\cite[p. 95]{Wash}.)
\begin{prop} Suppose $l\geq 3$ and for each, $r|l$, let $\cores^{\bbQ(\mu_l)}_{\bbQ(\mu_r)}:\bbC[G(r)]\rightarrow\bbC[G(l)]$ be the corestriction map
defined as in Remark~\ref{rem: fattening D}. Then
\beq\label{eq: equivariant func eq}
a_{\bbQ(\mu_l)/\bbQ, S_l}^{-,\ast}=\frac{1}{l}\sum_{r|l\atop r\neq 1}\cores^{\bbQ(\mu_l)}_{\bbQ(\mu_r)}(\cA_r\theta_{\bbQ(\mu_r)/\bbQ, S_r})
\eeq
where $\cA_r$ denotes the `equivariant Gauss Sum' ${\displaystyle\sum_{g\in G(r)}g(\xi_r)g=
 \sum_{a=1\atop (a,r)=1}^r} \xi_r^a\sigma_{a,r}\in \bbC[G(r)]$.
\end{prop}
\bPf\ We sketch two alternatives. The first uses the much more general equivariant functional equation coming from Thms.~2.2 and~2.1 of~\cite{twizo}:
In the notations of that paper, take $k=\bbQ$ and $\fm$ to be the cycle
$(l\bbZ)\infty$. Then, equations~(13) and~(9) of~{\em ibid.}\ with $s=1$ and $0$ respectively, show that
$a_{\bbQ(\mu_l)/\bbQ, S_l}^{-,\ast}=\frac{1}{l}\frac{1-c}{2}\Phi_\fm(0)^\ast$ (if $l=2\tilde{l}$ with $\tilde{l}$ odd, we need also~(8) of {\em ibid.}).
The reader may check that equation~(15) of~{\em ibid.}\ with $s=0$ then gives~(\ref{eq: equivariant func eq}) above.

Alternatively, and more directly, let $\chi$ be any odd character of $G(l)$
linearly extended to $\bbC[G(l)]$, let $\hat{\chi}$ be the associated, {\em primitive} Dirichlet character modulo $f_\chi$ (which
divides $l$) and let
$T_\chi$ denote the set or primes $q$ dividing $l$ but not $f_\chi$. One shows that if $r$
is of the form $f_\chi\prod_{q\in T}q$ for some $T\subset T_\chi$ then
\[
\chi\left(\cores^{\bbQ(\mu_l)}_{\bbQ(\mu_r)}(\cA_r)\right)=\frac{\varphi(l)}{\varphi(r)}\prod_{q\in T}(-\hat{\chi}(q))\tau(\hat{\chi})
\]
(where $\tau(\hat{\chi})$ is the usual Gauss Sum and $\varphi$ is Euler's function) and
otherwise $\chi(\cores^{\bbQ(\mu_l)}_{\bbQ(\mu_r)}(\cA_r))=0$.
Using this fact and some further manipulation, one can evaluate $\chi(\mbox{R.H.S. of~(\ref{eq: equivariant func eq})})$ in terms of
$L(0,\hat{\chi}\inv)$. The usual functional equation for $L(s,\hat{\chi})$ then shows
that it is precisely equal to $(i/\pi)L_{\bbQ(\mu_l)/\bbQ, S_l}(1,\chi)=\chi(\mbox{L.H.S. of~(\ref{eq: equivariant func eq})})$.
Since $\chi$ was an arbitrary odd character and both sides of~(\ref{eq: equivariant func eq}) lie in
$\bbC[G(l)]^-$, the result follows.\ePf
\noindent We now specialise to the case $F=K_n$ for $n\geq 0$ for our fixed but general 
abelian field $K$.
We shall always take $S$ to be $S^0(K_n/\bbQ)$ which equals $S_{f_n}\cup\{\infty\}$ and contains $p$. It is
independent of $n\geq 0$ so we drop it from the notation. It follows easily from the definition
that $a^-_{K_n/\bbQ}$ is the image of $a_{\bbQ(\mu_{f_n})/\bbQ}^-$ under the restriction map $\bbC[\Gal(\bbQ(\mu_{f_n})/\bbQ)]\rightarrow\bbC[G_n]$
coming from the inclusion $K_n\subset \bbQ(\mu_{f_n})$. In principle, a rather complicated formula for $a^-_{K_n/\bbQ}$ then follows
from~(\ref{eq: equivariant func eq}). A much simpler one -- which is also better suited to present purposes --
is easily obtained by the same process from a different formula for $a_{\bbQ(\mu_{f_n})/\bbQ}^-$ proved in~\cite[Lemma 7.1 (ii)]{laser}. (Note: $a_{\bbQ(\mu_{f_n})/\bbQ}^-$ would there be denoted
$a^-_{K_{f_n}/\bbQ,S}$.) The reader may check that this results in the following expression, which
shows in particular that $a^-_{K_n/\bbQ}$ lies in $K_n[G_n]^-$.
\beq\label{eq: def of aFminus}
a^-_{K_n/\bbQ}=\frac{1}{2f_n}(1-c)\sum_{g\in G_n}g\left(\Tr_{\bbQ(\mu_{f_n})/K_n}\left(\frac{\xi_{f_n}}{1-\xi_{f_n}}\right)\right)g\inv
\eeq
For each $\fP\in S_p(K_n)$ we shall write
$K_{n,\fP}$ for the (abstract) completion of $K_n$ at $\fP$. We shall usually regard the canonical embedding
$i_\fP:K_n\rightarrow K_{n,\fP}$ as an inclusion.
We write $K_{n,p}$ for the product $\prod_{\fP\in S_p(K_n)}K_{n,\fP}$ in
which we shall usually consider $K_n$ to be diagonally embedded ({\em via} $\prod_\fP i_\fP$).
Let $\pi_\fP$ denote the projection from $K_{n,p}$ to $K_{n,\fP}$ and  $U^1(K_{n,p})$ the group of `principal $p$-semilocal
units of $K_n$' \ie\ $\prod_{\fP\in S_p(K_n)}U^1(K_{n,\fP})$
considered as a multiplicative pro-$p$ group under the product topology.
({\em Warning:} we shall sometimes write it additively.)
$K_{n,p}$ is equipped with a natural $G_n$-action extending that on $K_n$
(see \eg~\cite[\S 2.3]{laser}) and such that
$U^1(K_{n,p})$ identifies as a finitely-generated, topological
$\bbZ_p[G_n]$-module with the
Sylow pro-$p$ subgroup of $(\cO_{K_n}\otimes_\bbZ \bbZ_p)^\times$.
We fix once and for all an algebraic closure $\barbbQ_p$ of
$\bbQ_p$ and an embedding $j:\barbbQ\rightarrow\barbbQ_p$ whose restriction to $K_n$ extends to an embedding
$j:K_{n,\fP^0}\rightarrow \barbbQ_p$ for some $\fP^0\in S_p(K_n)$. We shall also write $j$ for the composite $j\circ\pi_{\fP^0}$
taking $K_{n,p}$ onto $\overline{j(K_n)}$ (topological closure). We write
$\log_p$ for the $p$-adic logarithm defined by the usual convergent series on
$U^1(\overline{j(K_n)})$ and on $U^1(K_{n,\fP})$ for any $\fP$.
Given any $u\in U^1(K_{n,p})$ we set
$\lambda_{p,n}(u):=\sum_{g\in G_n}\log_p(j(gu))g\inv\in \overline{j(K_n)}[G_n]$.
Applying $j$ coefficientwise to $a^{-,\ast}_{K_n/\bbQ}$, we get an element $j(a^{-,\ast}_{K_n/\bbQ})$ of $j(K_n)[G_n]^-$ and a a map
\displaymapdef{\fs_n}{U^1(K_{n,p})}{\bbQ_p[G_n]^-}{u}{j(a^{-,\ast}_{K_n/\bbQ})\lambda_{p,n}(u)}
This is the map $\fs_{K_n/\bbQ,S^0(K_n/\bbQ)}$ of~\cite{laser} (taking
`$\tau_1$' to be $1\in\Gal(\barbbQ/\bbQ)$). The fact that $\fs_n(u)$
has coefficients in $\bbQ_p$ and is independent of $j$ therefore follows from~\cite[Prop.~2.16]{laser} (or, in our special case,
from~(\ref{eq: explicit formula for s(u)}) below).
It is clearly is $\bbZ_p[G_n]$-linear on $U^1(K_{n,p})$ and so
factors through the projection on $U^1(K_{n,p})^-$. Assuming $u$ lies in $U^1(K_{n,p})^-$, the
formula~(\ref{eq: def of aFminus}) gives
\beq\label{eq: explicit formula for s(u)}
\fs_n(u)=\sum_{g\in G_n}\sum_{\fP\in S_p(K_n)}\frac{1}{f_n}
{\rm Tr}_{K_{n,\fP}/\bbQ_p}\left(
\Tr_{\bbQ(\mu_{f_n})/K_n}\left(\frac{\xi_{f_n}}{1-\xi_{f_n}}\right)\log_p(\pi_\fP(g\inv u))
\right)g
\eeq
The next result gives the properties of $\fs_n$ that are crucial to the present paper. For
each $\fP\in S_p(K_n)$ we write
$(\cdot,\cdot)_{K_{n,\fP},\pnpo}$ for the Hilbert symbol on $K_{n,\fP}^\times\times K_{n,\fP}^\times$ with
values in $\mu_\pnpo$ (regarded as a subgroup of $K_{n,\fP}^\times$) defined as in~\cite{Neukirch}. This gives rise to a $\cR_n$-valued pairing
$[\cdot,\cdot]_{\fP,n}$ on $K_{n,\fP}^\times\times K_{n,\fP}^\times$ defined by
\[
\zeta_n^{[\alpha,\beta]_{\fP,n}}=(\alpha,\beta)_{K_{n,\fP},\pnpo}
\]
and hence, letting $\fP$ vary, to a pairing
\displaymapdef{[\cdot,\cdot]_n}{K_{n,p}^\times\times K_{n,p}^\times}{\cR_n}{(\alpha,\beta)}
{\displaystyle\sum_{\fP\in S_p(K_n)}[\pi_\fP(\alpha),\pi_\fP(\beta)]_{\fP,n}}
Properties of the Hilbert symbol give the following (see~\cite[eq. (18)]{laser})
\beq\label{eq: gal prop of Hilbert Symbol}
[g\alpha,g\beta]_n=\chi_{cyc,n}(g)[\alpha,\beta]_n \ \ \ \ \mbox{for all $\alpha,\beta\in K_{n,p}^\times$ and $g\in G_n$}
\eeq
Now write $\mu_{p^\infty}(K_{n,p})$ for $\prod_{\fP\in S_p(K_n)}\mu_{p^\infty}(K_{n,\fP})={\rm tor}_{\bbZ_p}(U^1(K_{n,p}))$
and $\fS_n$ for the image of $\fs_n$ in $\bbQ_p[G_n]^-$ (denoted $\fS_{K_n/k,S^0(K_n/\bbQ)}$ in~\cite{laser}).
\begin{prop}\label{prop: properties of sn} For all $n\geq 0$ we have
\begin{enumerate}
\item \label{part: properties of sn 1}
$
\ker (\fs_n|_{U^1(K_{n,p})^-})=\mu_{p^\infty}(K_{n,p})^-
$
\item \label{part: properties of sn 2}
 $\fS_n$ is contained in $\bbZ_p[G_n]^-$ with finite index and 
\beq\label{eq: CC for KnoverQ}
\fs_n(u)\equiv -{\textstyle \half}\sum_{g\in G_n}[\eps_n, g\inv u]_n g\ \ \ \ \mbox{modulo $\pnpo$, for all $u\in U^1(K_{n,p})^-$.}
\eeq
\end{enumerate}
\end{prop}
\bPf\ Part~\ref{part: properties of sn 1} follows easily from the fact $a^{-,\ast}_{K_n/\bbQ}$ is a unit of $\barbbQ[G_n]^-$
(since $\chi(a^{-,\ast}_{K_n/\bbQ})=(i/\pi)L_{K_n/\bbQ}(1,\chi)\neq 0$ for all odd $\chi\in \hat{G}_n$) or as a
special case of~\cite[Prop. 2.17]{laser} with $d=[k:\bbQ]=1$. The latter also shows $\bbQ_p\fS_n=\bbQ_p[G_n]^-$.  It
remains to show $\fS_n\subset\bbZ_p[G_n]^-$ and~(\ref{eq: CC for KnoverQ}). But it is easy to see that these amount precisely to the case of the
{\em `Congruence Conjecture'} of~\cite[\S 3]{laser} with
data $K_n/\bbQ,S=S^0(K_n/\bbQ)=S^1(K_n/\bbQ),p$ and $n$, which was proven in \ibid, Theorem~4.3. In particular,~(\ref{eq: CC for KnoverQ}) follows from
equations~(24) and~(20) of~\cite{laser}, taking $d=1$ and $\tau_1=1$ and noting that
`$\eta_{K_n^+/\bbQ,S^1(K_n/\bbQ)}$' equals our $-\half\otimes\eps_n$. (This last equation is established in the case $K_n=\bbQ(\mu_{f_n})$
during the course of the proof of \ibid, Theorem~4.3, see pp.~177-178. The general case follows on applying $N_{\bbQ(\mu_{f_n})^+/K_n^+}$ to both sides
and using~\cite[Prop. 5.7]{laser}.)\ePf
\rem\label{TEST} For those unfamiliar with~\cite{laser}, the following may shed some light on~(\ref{eq: CC for KnoverQ}).\\
\hspace*{3ex}(i)~The R.H.S.\ is $G_n$-equivariant in $u\in U^1(K_{n,p})$ and
lies in the {\em minus}-part of $\cR_n[G_n]$. (Use~(\ref{eq: gal prop of Hilbert Symbol}) with $g=c$).
Thus~(\ref{eq: CC for KnoverQ}) would read `$0\equiv 0$' for $u\in U^1(K_{n,p})^+$.\\
\hspace*{3ex}(ii)~If $K=\bbQ$ then $K_n=\bbQ(\mu_\pnpo)$ and $\eps_n=(1-\zeta_n)(1-\zeta_n\inv)$. In this case
the reader can easily check that~(\ref{eq: CC for KnoverQ}) follows immediately from~(\ref{eq: explicit formula for s(u)}) and the
explicit reciprocity law of Artin and Hasse,~\cite{Artin-Hasse}. Coleman's generalisation of
this law in~\cite{Coleman Dilogarithm} is an essential ingredient in the proof of Theorem~4.3 of~\cite{laser} which
establishes~(\ref{eq: CC for KnoverQ}) in the general case.\vertsp\\
For each $m\geq n\geq 0$, the norm $N^m_n:K_m^\times\rightarrow K_n^\times$ is the restriction of the
map $K_{m,p}^\times\rightarrow K_{n,p}^\times$ which is given by the products of local norms
(and also denoted $N^m_n$). Proposition 5.5 of~\cite{laser} gives a commuting diagram
\[
\xymatrix{
U^1(K_{m,p})\ar[dd]_{N^m_n}\ar[rrr]^{\fs_m}&&&\bbZ_p[G_m]^-\ar[dd]^{\pi^m_n}\\
&&&&&\\
U^1(K_{n,p})\ar[rrr]^{\fs_n}&&&\bbZ_p[G_n]^-\\
}
\]
We write $U^1_\infty$ for the projective limit of
the groups $U^1(K_{n,p})$ for all $n\geq 0$ with respect to the maps $N^m_n$, considered
as a natural $\Lambda_G$-module. The maps $(\fs_n)_{n\geq 0}$ give rise to a $\Lambda_G$-linear map
$\fs_\infty:U_\infty^1\rightarrow\Lambda_G^-$ factoring
through $U_\infty^{1,-}$.
The image of $\fs_\infty$ is precisely $\fS_\infty:={\displaystyle\lim_\leftarrow\fS_n}$ considered as a submodule of
$\Lambda_G^-$. (This follows from the finiteness of $\mu_{p^\infty}(K_{n,p})$ and Lemma~15.16 of~\cite{Wash} or the fact that
$N_n^m:\mu_{p^\infty}(K_{m,p})\rightarrow\mu_{p^\infty}(K_{n,p})$ is surjective
for all $m\geq n\geq 0$.) So, by Prop.~\ref{prop: properties of sn}~\ref{part: properties of sn 1} we obtain
an exact sequence of $\Lambda_G$-modules
\beq\label{eq: es for sinfty}
0\rightarrow \mu_{local,\infty}^-\hookrightarrow U^{1,-}_\infty \stackrel{\fs_\infty}{\longrightarrow}
\Lambda_G^-\longrightarrow \Lambda_G^-/\fS_\infty\rightarrow 0
\eeq
where $\mu_{local,\infty}$ denotes the projective limit of $\mu_{p^\infty}(K_{n,p})$  with respect to $N^m_n$ for all $m\geq n\geq 0$.

On the other hand, for each $\fP\in S_p(K_n)$ the reciprocity map of local class field theory restricts to
a map $\psi_{n,\fP}$ from $U^1(K_\fP)$ onto the inertia subgroup above $\fP$ in $\fX_n$ so that
the product $\prod_{\fP\in S_p(K_n)}\psi_{n,\fP}$ defines a $\bbZ_p[G_n]${\em -equivariant} map
$\psi_n:U^1(K_{n,p})\rightarrow\fX_n$ with image
$\Gal(M_n/L_n)$. Global class field theory
(and the fact that $K_n$ is CM) show that $\ker (\psi_n|_{U^1(K_{n,p})^-})=\mu_{p^\infty}(K_n)\subset\mu_{p^\infty}(K_{n,p})^-$ so
we get an exact sequence of $\bbZ_p[G_n]$-modules
\[
0\rightarrow \mu_{p^\infty}(K_n)\hookrightarrow U^1(K_{n,p})^- \stackrel{\psi_n}{\longrightarrow} \fX_n^-\longrightarrow X_n^-\rightarrow 0
\]
for each $n\geq 0$. Furthermore, if  $m\geq n$, one has $\rho_n^m\circ\psi_m=\psi_n\circ N^m_n$ so, on passing to limits,
we obtain a $\Lambda_G$-linear map
$\psi_\infty:U_\infty^1\rightarrow\fX_\infty$. It is easy to see
that $\psi_\infty(U^{1,-}_\infty)=\Gal(M_\infty/L_\infty)^-$ so we get an
exact sequence of $\Lambda_G$-modules
\beq\label{eq: es for psiinfty}
0\rightarrow \mu_{global,\infty}\hookrightarrow U^{1,-}_\infty \stackrel{\psi_\infty}{\longrightarrow} \fX_\infty^-\longrightarrow X_\infty^-\rightarrow 0
\eeq
where $\mu_{global,\infty}$ denotes the projective limit of $\mu_{p^\infty}(K_n)$ with respect to $N^m_n$ for all $m\geq n\geq 0$.
\begin{thm}\label{thm: commuting triangle with j, s and psi}
The following diagram commutes.
\[
\xymatrix{
U_\infty^1\ar[dd]_{\psi_\infty}\ar[drrrr]^{\fs_\infty}&&&&\\
&&&&\Lambda_G^-\\
\fX_\infty\ar[urrrr]_{\fj_\infty}&&&&
}
\]
\end{thm}
\bPf\
Suppose
$m\geq 0$ and let $v=(v_\fP)_{\fP\in S_p(K_m)}$ be an element of $U^1(K_{m,p})$. Then
\begin{eqnarray*}
\zeta_m^{\langle\bar{\eps}_m,\overline{\psi_m(v)}\rangle_m}&=&\psi_m(v)(\eps_m^{1/p^{m+1}})/\eps_m^{1/p^{m+1}}\\
&=&\prod_{\fP\in S_p(K_m)}\left(\psi_{m,\fP}(v_\fP)(\eps_m^{1/p^{m+1}})/\eps_m^{1/p^{m+1}}\right)\\
&=&\prod_{\fP\in S_p(K_m)}(v_\fP,\eps_m)_{K_{m,\fP},p^{m+1}}\\
&=&\prod_{\fP\in S_p(K_m)}(\eps_m,v_\fP)\inv_{K_{m,\fP},p^{m+1}}\\
&=&\zeta_m^{-[\eps_m,v]_m}
\end{eqnarray*}
(where the third equality comes from
the definition of the Hilbert symbol $(\cdot,\cdot)_{K_{m,\fP},p^{m+1}}$ that we are using and the fourth from
one of its basic properties).
Thus $\half\langle\bar{\eps}_m,\overline{\psi_m(v)}\rangle_m\equiv -\half[\eps_m,v]_m$ modulo $p^{m+1}$
and it follows from~(\ref{eq: def of equivariant pairing}) and~(\ref{eq: CC for KnoverQ}) that
\beq\label{eq: cong between j compose psi and s}
\{\bar{\eta}_m,\overline{\psi_m(u)}\}_m\equiv\fs_m(u) \pmod{p^{m+1}}\ \ \
\forall\,u\in U^1(K_{m,p})^-
\eeq
Hence if $\underline{u}=(u_m)_{m\geq 0}$ lies in $U^{1,-}_\infty$, we find
\[
\phi^m_n(\{\bar{\eta}_m,\overline{\psi_m(u_m)}\}_m)
\equiv\pi^m_n(\fs_m(u_m))\equiv\fs_n(u_n)\ \ \  \mbox{modulo $p^{m+1}$ for all $m\geq n\geq 0$.}
\]
Fixing $n$ and letting $m\rightarrow\infty$ gives
$\phi^\infty_n(\fj_\infty\circ\psi_\infty(\underline{u}))=\fs_n(u_n)=\phi^\infty_n(\fs_\infty(\underline{u}))$
in $\bbZ_p[G_n]^-$.
Since $n$ is arbitrary and both $\fj_\infty\circ\psi_\infty$ and
$\fs_\infty$ factor through $U^{1,-}_\infty$, the result follows.\ePf
\noindent Considering images and using $\fd_\infty=\iota_\infty\circ\fj_\infty$ and Theorem~\ref{thm: Stick in + part over Zp v2}~\ref{part: 1 of Stick in + part over Zp v2},
we deduce
\begin{cor}\label{cor: containments and annihilation}
\begin{enumerate}
\item \label{part: containments}
$\fJ_\infty\supset\fS_\infty$ and $\fD_\infty\supset\iota_\infty(\fS_\infty)$.
\item \label{part: annihilation}
$\iota_\infty(\fS_\infty)$ annihilates ${\displaystyle \lim_\leftarrow }\,A_m^+$ (or, equivalently, $X_\infty^+$).
 \ePf
\end{enumerate}
\end{cor}
(For comments on part~\ref{part: annihilation} in the case $K=\bbQ$, see Remark~\ref{rem: annihilation by fS vs Gras-Oriat}.)
Let $U^1(K_{n,p})^0:=\bigcap_{m\geq n}N^m_n(U^1(K_{m,p}))\subset U^1(K_{n,p})$. For any
$u\in U^1(K_{n,p})^-\cap U^1(K_{n,p})^0$ a compactness argument shows that we can find
$\underline{u}=(u_m)_m\in U^{1,-}_\infty$ with $u_n=u$.  Then
$\fs_n(u)=\phi^\infty_n(\fs_\infty(\underline{u}))=\phi^\infty_n(\fj_\infty(\psi_\infty(\underline{u})))$ by the Theorem. Using the
definition of $\fj_n$ we deduce
\begin{cor} $\fj_n(\psi_n(u))=\fs_n(u)$ for any $n\geq 0$ and $u\in U^1(K_{n,p})^-\cap U^1(K_{n,p})^0$.\ePf
\end{cor}
\noindent One might ask whether $\fj_n\circ\psi_n=\fs_n$ on the whole of $U^1(K_{n,p})^-$ and in particular whether $\fJ_n\supset\fS_n$.
(In general,
(\ref{eq: cong between j compose psi and s}) yields only the congruence $\fj_n\circ\psi_n\equiv\fs_n\pmod{p^{n+1}}$ on $U^1(K_{n,p})^-$.)
A sufficient (but possibly unnecessary) condition is that $U^1(K_{n,p})^-\subset U^1(K_{n,p})^0$. This clearly also
guarantees that $\phi^\infty_n(\fS_\infty)=\fS_n$.
\begin{lemma} Suppose $m> n\geq 0$. Then $N^m_n:U^1(K_{m,p})^-\rightarrow U^1(K_{n,p})^-$ is surjective iff $m\leq n_0$ or $c\in D_p(K_0/\bbQ)$.
\end{lemma}
\bPf\ let $T^m_n$ denote the (wild) inertia subgroup of $\Gal(K_m/K_n)$ at  primes of $K_n$ above $p$, on
which $D_n:=D_p(K_n/\bbQ)$-acts trivially.
Local class field theory gives an isomorphism of $\bbZ_p[G_n]$-modules between $U^1(K_{n,p})/N^m_nU^1(K_{m,p})$ and
$T^m_n\otimes_{\bbZ_p[D_n]}\bbZ_p[G_n]\cong T^m_n\otimes_{\bbZ_p}\bbZ_p[G_n/D_n]$ with $G_n$ acting via the second factor.
If $m\leq n_0$ then $T_n^m=\{0\}$. Otherwise  $T_n^m\neq\{0\}$ and $\bbZ_p[G_n/D_n]^-=\{0\}\Leftrightarrow c\in D_n\Leftrightarrow c\in D_0$ (since
$[K_n:K_0]$ is a power of $p\neq 2$).\ePf
\noindent Thus  $U^1(K_{n,p})^-\subset U^1(K_{n,p})^0$ for some $n\geq 0$ if and only if $c\in D_p(K_0/\bbQ)$ which implies in turn that
$U^1(K_{n,p})^-\subset U^1(K_{n,p})^0$ for all $n\geq 0$.
From the above arguments, we deduce:
\begin{cor}\label{cor: fJn contains fSn} Suppose $c\in D_p(K_0/\bbQ)$ \ie\ the primes of $K^+_0$ above $p$ do not split in $K_0$.
Then $\fS_n$ equals $\phi^\infty_n(\fS_\infty)$ and is contained in
$\fJ_n$. In particular, $\fJ_n$ is of finite index in $\bbZ_p[G_n]^-$.\ePf
\end{cor}
\noindent Next, passing to the quotient in the exact sequences~(\ref{eq: es for sinfty}) and~(\ref{eq: es for psiinfty}) gives
injective maps $U^{1,-}_\infty/\mu_{local,\infty}^-\rightarrow\Lambda_G^-$ and $U^{1,-}_\infty/\mu_{global,\infty}\rightarrow\fX_\infty^-$ which we denote
$\bar{\fs}_\infty$ and $\bar{\psi}_\infty$ respectively.
\begin{thm}\label{thm: commdiag} There is a commuting diagram of $\Lambda_G$-modules with exact rows and columns:
\beq\label{eq: big diag}
\xymatrix{
&0 \ar[d]&0\ar[d]&0\ar[d]\\
0\ar[r]&\mu_{local,\infty}^-/\mu_{global,\infty}\ar[d]\ar[r]&U^{1,-}_\infty/\mu_{global,\infty}\ar[d]^{\bar{\psi}_\infty}\ar[r]&
U^{1,-}_\infty/\mu_{local,\infty}^-\ar[d]^{\bar{\fs}_\infty}\ar[r]&0\ar[d]\\
0\ar[r]&\Gal(M_\infty/N_\infty^0)^-\ar[d]\ar[r]&\fX_\infty^-\ar[d]\ar[r]^{\fj_\infty}&
\Lambda_G^-\ar[d]\ar[r]&\Lambda_G^-/\fJ_\infty\ar@{=}[d]\ar[r]&0\\
0\ar[r]&\Gal(L_\infty/L_\infty\cap N_\infty^0)^-\ar[d]\ar[r]&X_\infty^-\ar[d]\ar[r]^{\fj'_\infty}&
\Lambda_G^-/\fS_\infty\ar[d]\ar[r]&\Lambda_G^-/\fJ_\infty\ar[d]\ar[r]&0\\
&0&0&0&0}
\eeq
\end{thm}
\bPf\ The exactness of the second and third columns follows from~(\ref{eq: es for psiinfty}) and~(\ref{eq: es for sinfty}) respectively.
The commutativity of the top, middle square is Theorem~\ref{thm: commuting triangle with j, s and psi}. There is therefore a unique map
$\fj'_\infty:X_\infty^-\rightarrow\Lambda_G^-/\fS_\infty$ making the bottom, middle square commute. The exactness of the top row is tautologous
and that of the middle row follows from
Thm.~\ref{thm: ker of jinf and dinf}.
A diagram chase then
shows that the isomorphism  $U^{1,-}_\infty/\mu_{global,\infty}\rightarrow \Gal(M_\infty/L_\infty)^-$ induced by $\bar{\psi}_\infty$
takes $\mu_{local,\infty}^-/\mu_{global,\infty}$
{\em onto} $\Gal(M_\infty/L_\infty)^-\cap\Gal(M_\infty/N_\infty^0)^-=\Gal(M_\infty/L_\infty N_\infty^0)^-$.
The rest follows easily.\ePf
\noindent Let $M_\infty^-$ be as in Section~\ref{sec: var K and kernels}. Similarly, let  $L_\infty^-$ denote the fixed field
of  $\Gal(L_\infty/K_\infty)^+$ acting on $L_\infty$, so that
$\Gal(M_\infty/L_\infty N_\infty^0)^-$ maps isomorphically onto $\Gal(M_\infty^-/L_\infty^- N_\infty^0)$.
The above proof then gives the following (implicit in~(\ref{eq: big diag})).
\begin{cor}\label{cor: first of comm diag.} $\bar{\psi}_\infty$ induces a $\Lambda_G$-isomorphism
$\mu_{local,\infty}^-/\mu_{global,\infty}\cong\Gal(M_\infty^-/L^-_\infty N_\infty^0)$. \ePf
\end{cor}

We now give an explicit description of $\mu_{local,\infty}^-/\mu_{global,\infty}$ as a $\Lambda_G$-module.
Recall that $K_{n_0}=K_{m_0}=F_{m_0}$ for an integer $m_0\geq n_0$ and an abelian field $F$ unramified over $\bbQ$ above $p$.
(See Remark~\ref{rem: m_0, i_0, F etc.}.) Suppose that $n\geq m_0$. It follows that $K_n=F_n$, $\mu_{p^\infty}(K_n)=\mu_{\pnpo}$ and also
$\mu_{p^\infty}(K_{n,\fP})=i_\fP(\mu_\pnpo)$ for all $\fP\in S_p(K_n)$. Hence
\displaymapdef{\nu_n}{\bbZ_p[S_p(K_n)]\otimes_{\bbZ_p}\mu_\pnpo}{\mu_{p^\infty}(K_{n,p})}
{\displaystyle\sum_{\fP\in S_p(K_n)}a_\fP\fP\otimes\zeta_\fP}{(i_\fP(\zeta_\fP)^{a_\fP})_\fP}
is an isomorphism of $\bbZ_p[G_n]$-modules (where $g(\sum_{\fP}a_\fP\fP\otimes\zeta_\fP)=\sum_{\fP}a_\fP g(\fP)\otimes g(\zeta_\fP)$ for all $g\in G_n$).
Note that $\mu_{p^\infty}(K_n)$ is the image of $\bbZ_p(\sum_\fP\fP)\otimes\mu_\pnpo$ under $\nu_n$ and $\mu_{p^\infty}(K_{n,p})^-$ is that
of $(\bbZ_p[S_p(K_n)]\otimes\mu_\pnpo)^-=((1+c)\bbZ_p[S_p(K_n)])\otimes\mu_\pnpo$.
Since $K_\infty/K_{n_0}$ is totally ramified above $p$, we can identify $\bbZ_p[S_p(K_n)]$ with $\bbZ_p[S_p(K_{n_0})]$ and, hence,
$(1+c)\bbZ_p[S_p(K_n)]$ with $\bbZ_p[S_p(K_{n_0}^+)]$ for any $n\geq m_0$. Moreover, if $m\geq n\geq m_0$ then $N_n^m: \mu_{p^\infty}(K_{m,p})
\rightarrow\mu_{p^\infty}(K_{n,p})$  is simply the $p^{m-n}$th power map. Passing to the limit and then the quotient we find easily
\beq\label{eq: double iso}
\left(\frac{\bbZ_p[S_p(K_{n_0}^+)]}{\bbZ_p(\sum_\fP\fP)}\right)\otimes_{\bbZ_p} \bbZ_p(1)
\cong\mu_{local,\infty}^-/\mu_{global,\infty}
\cong\Gal(M_\infty^-/L_\infty^- N_\infty^0)
\eeq
as $\Lambda_G$-modules, where $\bbZ_p(1)$ is a rank-1 $\bbZ_p$-module with $G_\infty$ acting through $\chi_{cyc}$.
Using also Corollary~\ref{cor: what happens when Sp has one element I} we deduce:
\begin{cor}\label{cor: what happens when Sp has one element II}
\begin{enumerate}
\item $|S_p(K_{n_0}^+)|=1\Leftrightarrow M_\infty^-=L_\infty^-N_\infty^0$
\item Suppose $|S_p(K_{n_0}^+)|=1$. Then $N_\infty^0=N_\infty$, $M_\infty^-=L_\infty^-N_\infty$ and
\beq
(\Gal(L_\infty/L_\infty\cap N_\infty)^-)^\dag\cong
(\Gal(M_\infty/N_\infty)^-)^\dag\cong
 \Hom_{\bbZ_p}\left(A^+_\infty,\bbQ_p/\bbZ_p\right)
\eeq
as $\Lambda_G^+$-modules. In particular, Greenberg's Conjecture holds for $K_\infty^+/K^+$ if and only if $L_\infty^-\subset N_\infty$.\ePf
\end{enumerate}
\end{cor}
\noindent Note that the equality $M_\infty^-=L_\infty^-N^-_\infty$
appears to be known already in certain cases, even without the condition $|S_p(K_{n_0}^+)|=1$.

\noindent We write $\Lambda_\Gamma$ for $\bbZ_p[[\Gamma_0]]$ which is isomorphic to $\bbZ_p[[\bbZ_p]]\cong\bbZ_p[[X]]$ and so is a noetherian
integral domain (often denoted $\Lambda$). Clearly, $\Lambda_G$ is a $\Lambda_\Gamma$-algebra and it is easy to see that
$\Lambda_G^-$ is free
of rank $\half[K_0:\bbQ]$ over $\Lambda_\Gamma$. We may consider~(\ref{eq: big diag}) and~(\ref{eq: double iso}) over
$\Lambda_\Gamma$ by restriction of scalars.
Using also some `classical' results from  Iwasawa Theory, this yields:
\begin{thm}\label{thm: Lambda_Gamma Structure}
\begin{enumerate}
\item\label{part: 1 of Lambda_Gamma Structure}
All the modules in diagram~(\ref{eq: big diag}) are finitely generated over $\Lambda_\Gamma$.

\item\label{part: 1.5 of Lambda_Gamma Structure}

Those in the left-hand column and the bottom row
are $\Lambda_\Gamma$-torsion.

\item\label{part: 2 of Lambda_Gamma Structure}
$\Gal(M_\infty/N_\infty^0)=\ker(\fj_\infty)$ is precisely ${\rm tor}_{\Lambda_\Gamma}(\fX_\infty)$ (the $\LG$-torsion submodule of $\fX_\infty$).
\end{enumerate}
\end{thm}
\bPf\ Recall that a $\LG$-module $A$ is said to be {\em pseudo-isomorphic} to another, $B$ (written $A\sim B$)
iff there exists a $\LG$-homomorphism $A\rightarrow B$ with finite
kernel and cokernel. It is shown in~\cite{Wash} that $\fX_\infty$ is finitely generated over $\LG$ (see~p.~292-293) and that
\[
\fX_\infty\sim \LG^{\half[K_0:\bbQ]}\oplus C
\]
for some finitely generated torsion $\LG$-module $C$ (see Theorem~13.31 \ibid). Thus both $\fX_\infty^-$ and $\Lambda_G^-$ are
finitely generated and~\ref{part: 1 of Lambda_Gamma Structure} follows.
Next, it is well-known that $X_\infty$ is finitely generated
and torsion over $\LG$ (see \eg\ \cite[\S 13.3]{Wash}). Furthermore, we have
\beq\label{eq: second for  Lambda_Gamma Structure}
X_\infty\sim \Hom_{\bbZ_p}(A_\infty,\bbQ_p/\bbZ_p)\cong\Gal(M_\infty/N_\infty)^\dag
\eeq
as $\Lambda_\Gamma$-modules, where the isomorphism follows from our~(\ref{eq: gal-equiv pairing}) and the pseudo-isomorphism
from~\cite[Prop. 15.34]{Wash} ({\em Note:} Washington's `$\tilde{X}$' instead of `$X$' comes about because
because of his different action on $\Hom$'s.)  Since $N_\infty$ is contained in $M_\infty^-$ we have
$\fX_\infty^+\subset\Gal(M_\infty/N_\infty)$ and it follows from~(\ref{eq: second for  Lambda_Gamma Structure})
that $\fX_\infty^+$ is also $\Lambda_\Gamma$-torsion. Finally, eq.~(\ref{eq: double iso})
shows that $\mu_{local,\infty}^-/\mu_{global,\infty}$ is killed by $(\gamma_{n_0}-\chi_{cyc}(\gamma_{n_0}))\in\LG$.
It follows from the above facts that $X_\infty^-$ and all the modules in the L.H.\ column of~(\ref{eq: big diag}) are $\LG$-torsion
so part~\ref{part: 1.5 of Lambda_Gamma Structure} will follow if we can show that
$\Lambda_G^-/\fJ_\infty$ is too, \ie\ that $(\Lambda_G^-/\fJ_\infty)\otimes_{\LG}\cF_\Gamma=\{0\}$ where $\cF_\Gamma$
denotes the field of fractions of $\LG$. Consider the exact sequence obtained by applying $\otimes_{\LG}\cF_\Gamma$
to the middle row of of~(\ref{eq: big diag}). From the torsion results proved so far, the first term of this sequence vanishes and so does
$\fX_\infty^+\otimes_{\LG}\cF_\Gamma$. Thus $\fj_\infty\otimes 1$ is injective and
\begin{eqnarray*}
\dim_{\cF_\Gamma}( \fX_\infty^-\otimes_{\LG}\cF_\Gamma)&=&
   \dim_{\cF_\Gamma}( \fX_\infty^-\otimes_{\LG}\cF_\Gamma)+\dim_{\cF_\Gamma}(\fX_\infty^+\otimes_{\LG}\cF_\Gamma)\\
&=&\dim_{\cF_\Gamma}( \fX_\infty\otimes_{\LG}\cF_\Gamma)\\
&=&\half[K_0:\bbQ]=\dim_{\cF_\Gamma}(\Lambda_G^- \otimes_{\LG}\cF_\Gamma)
\end{eqnarray*}
Hence $\fj_\infty\otimes 1$ is also surjective and~\ref{part: 1.5 of Lambda_Gamma Structure} follows.
For part~\ref{part: 2 of Lambda_Gamma Structure}, we already know that $\Gal(M_\infty/N_\infty^0)^+=\fX_\infty^+$ and
$\Gal(M_\infty/N_\infty^0)^-$ are $\LG$-torsion, so $\ker(\fj_\infty)=\Gal(M_\infty/N_\infty^0)\subset {\rm tor}_{\Lambda_\Gamma}(\fX_\infty)$.
The reverse inclusion, ${\rm tor}_{\Lambda_\Gamma}(\fX_\infty)\subset\ker(\fj_\infty)$, is clear, since
$\Lambda_G^-$ is $\LG$-torsionfree.\ePf

\rem\ The equality $\Gal(M_\infty/N_\infty^0)={\rm tor}_{\Lambda_\Gamma}(\fX_\infty)$
appears to be known already if $f_{-1}$ is a power of $p$ (so $K_\infty=\bbQ(\mu_{p^\infty})$ and $N_\infty^0=N_\infty$)
but to be new for general abelian $K$. It's also worth pointing out that since $\Lambda_G$ is $\LG$-torsionfree and $\ker(\fj_\infty)=\Gal(M_\infty/N_\infty^0)$ is
$\Lambda_\Gamma$-torsion, the latter must be  the full right kernel of the pairing $\{\ ,\ \}_\infty$. In other words, if $h\in\fX_\infty$ then
$\{\underline{\alpha} ,h\}_\infty=0$ holds for all $\underline{\alpha}\in \cV_\infty$ iff $\{\underline{\eta} ,h\}_\infty=0$.
\vertsp\\
\noindent
\rem\ Suppose for a moment that $K$ is a general number field, not necessarily abelian over $\bbQ$. We can still define $N_\infty$, $N^0_\infty$,
$\Lambda_\Gamma$ \etc\ in much the same way, and also a subfield $T_\infty$ of $M_\infty$ by
$\Gal(M_\infty/T_\infty)={\rm tor}_{\Lambda_\Gamma}(\fX_\infty)$. The generalisation
of~(\ref{eq: second for  Lambda_Gamma Structure}) implies
$T_\infty\subset N_\infty$ and it would be interesting to know whether
one has $T_\infty=N^0_\infty$, as in the abelian case.

In the general case, the fixed field of
$\bar{\psi}_\infty(\mu_{local,\infty}/\mu_{global,\infty})$ acting on $M_\infty$
is the {\em field  of Bertrandias-Payan} denoted $K^{BP}_\infty$ (see \eg~\cite
{NQD's second torsion paper}). For $K$ abelian, Cor.~\ref{cor: first of comm diag.} and
Theorem~\ref{thm: Lambda_Gamma Structure}~
\ref{part: 2 of Lambda_Gamma Structure} give respectively
\[
K^{BP}_\infty\cap M_\infty^-=L_\infty^-N_\infty^0=L_\infty^-T_\infty
\]
For $K$ general (but CM) the equality of the first and last terms is clearly equivalent to the
$\Lambda_\Gamma$-torsionfreeness of $U^{1,-}_\infty/\mu_{local,\infty}^-$. The latter follows (at least in certain cases) from
work of Coleman. If $K$ is abelian, then of course it is a consequence of the injectivity of $\bar{\fs}_\infty$ in~(\ref{eq: big diag})
which was also crucial to Cor.~\ref{cor: first of comm diag.} \etc\ We stress that this injectivity was in turn deduced
from Prop.~\ref{prop: properties of sn}~\ref{part: properties of sn 1} and hence, at base, from the non-vanishing of {\em complex}
$L$-functions at $s=1$.
\vertsp\\
\noindent
Finally, we can apply $\underline{\ }^\dag$ to the bottom row of~(\ref{eq: big diag}). Since $\iota_\infty$
induces $\Lambda_G$-isomorphisms
$(\Lambda_G^-/\fS_\infty)^\dag\rightarrow\Lambda_G^+/\iota_\infty(\fS_\infty)$ and
$(\Lambda_G^-/\fJ_\infty)^\dag\rightarrow\Lambda_G^+/\fD_\infty$, we deduce:
\begin{cor} There is an exact sequence of $\Lambda_G^+$-modules which are f.g.\ and torsion over $\LG$:
\beq\label{eq: 4-term exact seq}
\xymatrix{
0\ar[r]&(\Gal(L_\infty/L_\infty\cap N_\infty^0)^-)^{\dag}\ar[r]&(X_\infty^-)^\dag\ar[rr]^(0.45){\iota_\infty\circ\fj'_\infty}&&
\Lambda_G^+/\iota_\infty(\fS_\infty)\ar[r]&\Lambda_G^+/\fD_\infty\ar[r]&0\\
}
\eeq
\ePf
\end{cor}
In the next section we shall analyse this
 sequence under the assumption $K=\bbQ$ which
eliminates at a stroke many of the complicating (but also interesting) phenomena of the general case.
The weaker assumption $p\ndiv[K:\bbQ]$ would allow us to decompose~(\ref{eq: big diag}) and~(\ref{eq: 4-term exact seq})
using ($p$-adic) characters of $G_0$ and hence `isolate' these
phenomena -- \eg\ the possible non-triviality of $\mu_{local,\infty}^-/\mu_{global,\infty}$ and/or
$\Gal(N_\infty/N_\infty^0)$ -- at certain `troublesome' characters.

\section{The Case $K=\bbQ$: Computation of $\fS_n$ and $\fS_\infty$, and the Main Conjecture}\label{sec: analysis of exact sequence}
We return to the situation and notations of Section~\ref{sec: d and D for K=Q}, so that $K_n=\bbQ(\mu_\pnpo)$.
We start by determining $\fS_n$ for $n\geq 0$ showing that in this case it is exactly the $\bbZ_p$-span of the Stickelberger ideal.
Firstly, $K_{n,p}$ is simply the completion of $K_n$ at its unique prime above $p$, so that
the embedding $j$ induces an isomorphism from $K_{n,p}$ to $\hat{K}_n:=\bbQ_p(j(\zeta_n))$.
For convenience we regard this as an identification and suppress $j$ from the notation. Thus,
$G_n=D_p(K_n/\bbQ)$
identifies with $\Gal(\hat{K}_n/\bbQ_p)$
whose action commutes with $\log_p$, giving
\beq\label{eq: new formula for sn}
\fs_n(u)=a^{-,\ast}_{K_n/\bbQ}\sum_{g\in G_n}g(\log_p(u))g\inv \ \ \ \ \mbox{for all $u\in U^1(\hat{K}_n)$}
\eeq
Let us write simply $\theta_n$ for the element
$\theta_{\bbQ(\mu_\pnpo)/\bbQ, S_{\pnpo}}$ of $\bbQ[G_n]^-$, $b_n$ for the
element $\frac{1}{\pnpo}\sum_{i=0}^n\zeta_n$ of $K_n\subset \hat{K}_n$
and $\bT_n$ for the trace pairing $\hat{K}_n\times
\hat{K}_n\rightarrow \bbQ_p$ that is, $\bT_n(v,w):=\Tr_{\hat{K}_n/\bbQ_p}(vw)\ \forall\, v,w\in \hat{K}_n$. $\bT_n$ is symmetric, non-degenerate
and clearly satisfies
\beq\label{eq: prop of pairing Tn}
\bT_n(xv,yw)=\bT_n(y^\ast xv,w)=\bT_n(v,x^\ast yw)\ \ \ \ \ \mbox{for all $v,w\in \hat{K}_n$ and $x,y\in \bbQ_p[G_n]$}
\eeq
We define a $\bbZ_p[G_n]$-equivariant map $\fw_n$ by
\displaymapdef{\fw_n}{U^1(\hat{K}_n)}{\bbQ_p[G_n]}{u}{\displaystyle \sum_{g\in G_n}\bT_n(b_n,g(\log_p(u)))g\inv}
\begin{prop}\label{prop: sn in terms of theta and w} $\fs_n(u)=\theta_n\fw_n(u)$ for all $u\in U^1(\hat{K}_n)$.
\end{prop}
\bPf The first equality in~(\ref{eq: def of theta})
(together with the fact $S_{p^{i+1}}=S_p\ \forall i$) implies
\beq\label{eq: pi applied to theta}
\pi^i_j(\theta_i)=\theta_j\ \ \ \ \mbox{For each $i>j\geq 0$}
\eeq
It follows easily from this, equation~(\ref{eq: equivariant func eq}) with $l=\pnpo$ and the definition of $\cA_{p^{i+1}}$ that
\begin{eqnarray*}
a_{K_n/\bbQ}^{-,\ast}&=&\frac{1}{\pnpo}\sum_{i=0}^n\cores^{K_n}_{K_i}(\cA_{p^{i+1}}\theta_i)\\
                     &=&\theta_n\frac{1}{\pnpo}\sum_{i=0}^n\cores^{K_n}_{K_i}(\cA_{p^{i+1}})\\
                     &=&\theta_n\sum_{h\in G_n}h(b_n)h
\end{eqnarray*}
Substituting this in~(\ref{eq: new formula for sn}) and rearranging gives the result.\ePf
\noindent The determination of the image of $\fw_n$ is a formal consequence of a
`classical' result from~\cite{Iwasawa on some mods...}: Let $\cL_n$ denote $\log_p(U^1(\hat{K}_n))$. This is easily seen to be
a $\bbZ_p[G_n]$-submodule of $\hat{K}_n$ of $\bbZ_p$-rank equal to $[\hat{K}_n:\bbQ_p]=|G_n|$.
Let $\cL_n^\star$ denote the $\bbZ_p$-dual of $\cL_n$ w.r.t.\ $\bT_n$, namely the set
$\{v\in \hat{K}_n: \bT_n(v,w)\in\bbZ_p\ \forall\,w\in \cL\}$. To determine $\cL_n^\star$ (which he denotes `$\fX_n$') Iwasawa defines a
fractional ideal $\fA_n$ of $\bbQ_p[G_n]$ which,
in our notation, is given by
\beq\label{eq: def of An}
\fA_n:=\bbZ_p[G_n]\left(-\theta_n^\ast+\frac{2}{p^n}\sum_{g\in G_n}g\right)+I(\bbZ_p[G_n])
\eeq
(See~\cite[p.~44]{Iwasawa on some mods...}. The element in large parentheses coincides with that denoted `$\xi_n$' by Iwasawa.) He also shows that
there is a $\bbQ_p[G_n]$-isomorphism $\bbQ_p[G_n]\rightarrow \hat{K}_n$  which he denotes `$\varphi_n$' and
which sends $x$ to $xcb_n$ in our notation. (See the start of \S 1.6, {\em ibid.}, noting that Iwasawa's
`$\theta_n$' is our $c(b_n)$.) Theorem~1 of~\cite{Iwasawa on some mods...} thus amounts in our notation to the equation
\beq\label{eq: Iwa's char. of dual of log}
\cL_n^\star=\fA_nb_n
\eeq
For any fractional ideal $\fC$ of $\bbQ_p[G_n]$, we define another fractional ideal $\fC^\star$ by
\[
\fC^\star:=\{y\in \bbQ_p[G_n]:x^\ast y\in \bbZ_p[G_n]\ \forall\,x\in \fC\}
\]
The reason for the similarity of notation is that $\fC^\star$ is easily seen to be the $\bbZ_p$-dual of $\fC$ w.r.t.\ the symmetric, non-degenerate
pairing
$\bB_n:\bbQ_p[G_n]\times\bbQ_p[G_n]\rightarrow\bbQ_p$ taking $(\sum_{g\in G_n}a_g g,\sum_{g\in G_n}b_g g)$ to
$\sum_{g\in G_n}a_g b_g$, \ie\  the coefficient of $1$ in $\left(\sum_{g\in G_n}a_g g\right)^\ast\left(\sum_{g\in G_n}b_g g\right)$. We can now prove
\begin{prop}\label{prop: image of w} $\im (\fw_n)=\fA_n^\star$.
\end{prop}
\bPf\ Eq.~(\ref{eq: Iwa's char. of dual of log}) shows that $\{gb_n:g\in G_n\}$ is a $\bbQ_p$-basis of $\hat{K}_n$ and it follows
from~(\ref{eq: prop of pairing Tn}) that the dual
basis w.r.t.\ $\bT_n$ is of form $\{hb'_n:h\in G_n\}$ for some $b'_n\in \hat{K}_n$. More precisely $\bT_n(gb_n,hb'_n)=\delta_{g,h}$ so that
$\bT_n(xb_n,yb'_n)=\bB_n(x,y)\ \forall\, x,y\in \bbQ_p[G_n]$. Now, clearly, $\cL_n$ must be of form
$\fC b'_n$ for some fractional ideal $\fC$ of $\bbQ_p[G_n]$. Since also
$\cL_n$ is the $\bbZ_p$-dual of $\cL_n^\star$ w.r.t.\ $\bT_n$, equation~(\ref{eq: Iwa's char. of dual of log}) gives,
for any $y\in \bbQ_p[G_n]$,
\[
y\in\fC\Leftrightarrow yb'_n\in\cL_n\Leftrightarrow \bT_n(xb_n,yb'_n)\in \bbZ_p\ \forall\, x\in \fA_n
\Leftrightarrow \bB_n(x,y)\in \bbZ_p\ \forall\, x\in \fA_n\Leftrightarrow y\in\fA_n^\star
\]
Thus $\cL_n=\fA_n^\star b'_n$. Finally, the map $\alpha:\hat{K}_n\rightarrow\bbQ_p[G_n]$ sending $v$ to $\sum_{g\in G_n}\bT_n(b_n,g(v))g\inv$
is clearly $\bbQ_p[G_n]$-equivariant, so $\im(\fw_n)=\alpha(\cL_n)=\alpha(\fA_n^\star b'_n)=\fA_n^\star\alpha(b'_n)=\fA_n^\star.1=\fA_n^\star$.\ePf
\noindent  Propositions~\ref{prop: sn in terms of theta and w}
and~\ref{prop: image of w} imply $\fS_n=\im(\fs_n)=\theta_n\im(\fw_n)=\theta_n\fA_n^\star$. Now
definition~(\ref{eq: def of An}) shows that $\fA_n$ and $\bbZ_p[G_n]\theta_n^\ast+\bbZ_p[G_n]$ have the same minus parts, hence so do
$\fA_n^\star$ and $(\bbZ_p[G_n]\theta_n^\ast+\bbZ_p[G_n])^\star=(\bbZ_p[G_n]\theta_n^\ast)^\star\cap\bbZ_p[G_n]^\star$. Since
also $\theta_n$ lies in $\bbQ_p[G_n]^-$, we deduce
\beq\label{eq: formula for fSn}
\fS_n
=\theta_n\left(\bbZ_p[G_n]\theta_n^\ast+\bbZ_p[G_n]\right)^\star
=\theta_n\{y\in\bbZ_p[G_n]:\theta_n y\in\bbZ_p[G_n]\}
\eeq
(Incidentally, this proves $\fS_n\subset\bbZ_p[G_n]^-$ independently of Prop~\ref{prop: properties of sn}.)
Since $\chi_{cyc}:G_\infty\rightarrow\bbZ_p^\times$ is an isomorphism, $G_\infty$ is pro-cyclic and we fix henceforth
a topological generator $g_{\infty}$ whose image $g_{n}$ in $G_n$ generates the latter.
\begin{lemma}
\beq\label{eq: equality of three ideals}
\bbZ_p[G_n](g_n-\chi_{cyc}(g_\infty))=\langle\sigma_{a,\pnpo}-a: (a,2p)=1\rangle_{\bbZ_p[G_n]}=
\{y\in\bbZ_p[G_n]:\theta_n y\in\bbZ_p[G_n]\}
\eeq
\end{lemma}
\bPf (Sketch). Denote the sets by $(1)$, $(2)$ and $(3)$ respectively. One checks directly that
$(2)\subset(3)$ and that $(3)/(2)$ is represented by elements of $(3)$ lying in $\bbZ_p$, which must clearly be divisible by $\pnpo$. Since
$-2\pnpo=\sigma_{1+2\pnpo,\pnpo}-(1+2\pnpo)$, we deduce $\pnpo\in (2)$ so $(3)=(2)$. Clearly,
$(1)$ is generated over $\bbZ_p[G_n]$ by
the elements $g_n^l-\chi_{cyc}(g_\infty^l)$ for $l\geq 1$. Taking $l=(p-1)p^n$ we find easily $\pnpo\in(1)$ so it
suffices to show $(1)\equiv(2)$ mod $\pnpo$.
But $g_n^l=\sigma_{a,\pnpo}$ implies $\chi_{cyc}(g_\infty^l)\equiv a\bmod \pnpo$ so the generators are the same mod $\pnpo$.\ePf
\rem\ In fact, if $G(r)$ denotes $\Gal(\bbQ(\mu_r)/\bbQ)$ for some $r>1$, it is well known that
\[
\ann_{\bbZ[G(r)]}(\mu(\bbQ(\mu_r)))=\langle\sigma_{a,r}-a: (a,2r)=1\rangle_{\bbZ[G(r)]}=
\{y\in\bbZ[G(r)]:\theta_{\bbQ(\mu_r)/\bbQ, S_r} y\in\bbZ[G(r)]\}
\]
The argument for the second equality above is similar to that for
the second equality
of~(\ref{eq: equality of three ideals}). For more details, see~\cite[Lemma 6.9]{Wash} but note
that the  element `$\theta$' there is our $-\theta_{\bbQ(\mu_r)/\bbQ, S_r}+\frac{1}{2}\sum_{g\in G(r)}g$.
The Stickelberger ideal ${\rm St}_{\bbQ(\mu_r)}$ of $\bbZ[G(r)]$ is
$\theta_{\bbQ(\mu_r)/\bbQ, S_r}\ann_{\bbZ[G(r)]}(\mu(\bbQ(\mu_r)))$. (This is the `unenlarged' ideal, but for $r=\pnpo$ it makes no difference.)
Thus eq.~(\ref{eq: formula for fSn}),  the second equality in~(\ref{eq: equality of three ideals}) and the first equality in the last equation
imply $\fS_n=\bbZ_p{\rm St}_{\bbQ(\mu_\pnpo)}$.\vertsp\\
Let $\tilde{\theta}_n=(g_n-\chi_{cyc}(g_\infty))\theta_n$. Then $\tilde{\theta}_n\in \bbZ_p[G_n]^-$ by~(\ref{eq: equality of three ideals}) and
the sequence $(\tilde{\theta}_n)_n$ defines an element $\tilde{\theta}_\infty$ of $\Lambda_G^-$ by~(\ref{eq: pi applied to theta}).
Equations~(\ref{eq: equality of three ideals}) and~(\ref{eq: formula for fSn})
give
\begin{thm}\label{prop: generation of fS} If $K=\bbQ$ then $\fS_n$ (for any $n\geq 0$) and $\fS_\infty$ are the principal
ideals of $\bbZ_p[G_n]^-$ and $\Lambda_G^-$ generated by
$\tilde{\theta}_n$ and $\tilde{\theta}_\infty$ respectively.\ePf
\end{thm}
\rem\ It is worth noting that a similarly simple description of $\fS_n$ cannot be expected
for general abelian $K$. Indeed, if $\theta_{K_n}$ denotes the
Stickelberger element of $\bbQ_p[G_n]^-$ generalising $\theta_n$ then the  phenomenon
of `trivial zeroes' means that $\bbZ_p[G_n]^-\cap\bbZ_p[G_n]^-\theta_{K_n}$ is frequently of infinite index in $\bbZ_p[G_n]^-$ and so cannot contain
$\fS_n$ which is always of finite index.\vertsp\\
\noindent
Our assumption $K=\bbQ$ implies $|S_p(K_0^+)|=1$ so, using Corollary~\ref{cor: what happens when Sp has one element II} and
Theorem~\ref{prop: generation of fS}, the sequence~(\ref{eq: 4-term exact seq}) can be rewritten
as 
\beq\label{eq: 4-term exact seq in case K=Q}
0\rightarrow A_\infty^{+,\vee}\longrightarrow(X_\infty^-)^\dag\longrightarrow
\Lambda_G^+/(\iota_\infty(\tilde{\theta}_\infty))\longrightarrow \Lambda_G^+/\fD_\infty
\rightarrow 0
\eeq
Since $p\ndiv |G_0|=p-1$, there is a is a unique
splitting $G_\infty=G_0\times \Gamma_0$ and we can
decompose~(\ref{eq: 4-term exact seq in case K=Q})
using (even) characters of $G_0$. 
Let $\omega: G_0\rightarrow\bbZ_p^\times$ be the Teichm\"uller character (the restriction of $\chi_{cyc}$) and let $e_{j}$ be the
idempotent of $\bbZ_p[G_0]$ associated to $\omega^j$ for $j\in\bbZ$. Any $\bbZ_p[G_0]$-module
$M$ is the direct sum of its components $M^{(j)}:=e_{j}M$ for $j=0,\ldots,p-2$.
For $\bbZ_p[G_0^+]$-modules, we restrict to $j$ even.
It follows that~(\ref{eq: 4-term exact seq in case K=Q})
is the direct sum of the exact sequences:
\beq\label{eq: 4-term exact seq in case K=Q with j}
0\rightarrow A_\infty^{(j),\vee}\longrightarrow(X_\infty^{(1-j)})^\dag\longrightarrow
\Lambda_G^{(j)}/(\iota_\infty(\tilde{\theta}_\infty))^{(j)}\longrightarrow
\Lambda_G^{(j)}/\fD_\infty^{(j)}
\rightarrow 0
\eeq
of f.g.\ torsion $\Lambda_\Gamma$-modules for $j=0,2,4,\ldots,p-3$.
The fact that the generalised Bernoulli number $B_{1,\omega\inv}$ lies in $p\inv\bbZ_p^\times$ implies that
the image of $e_0\iota_\infty(\tilde{\theta}_\infty)$ in $e_0\bbZ_p[G_0]=\bbZ_pe_0$
lies in $\bbZ_p^\times e_0$. It follows easily that $\Lambda_G^{(0)}/(\iota_\infty(\tilde{\theta}_\infty))^{(0)}$ vanishes;
but the same fact also
implies that $A_0^{(1)}\cong(X_\infty^{(1)})_{\Gamma_0}$ vanishes (by Stickelberger's Theorem) hence so does $X_\infty^{(1)}$.
Thus~(\ref{eq: 4-term exact seq in case K=Q with j}) is trivial for $j=0$ and we suppose henceforth $j\neq 0$ unless otherwise stated.
To analyse the
third non-zero term in~(\ref{eq: 4-term exact seq in case K=Q with j})
we first write
$g_\infty=g_0\gamma$ so that $\gamma$ and $\kappa:=\chi_{cyc}(\gamma)$ topologically
generate $\Gamma_0$ and $1+p\bbZ_p$ respectively.
Similarly $g_n=g_0\gamma(n)$ where $\gamma(n)$ is the
image of $\gamma$ in $\Gamma(n):=\Gal(K_n/K_0)\cong \Gamma_0/\Gamma_n$ and $G_n=G_0\times \Gamma(n)$.
Define $\tilde{\theta}_{n,j},\theta_{n,j}$ and $v_{n,j}\in \bbQ_p[\Gamma(n)]$ by
\[\mbox{$e_{1-j}\tilde{\theta}_n=\tilde{\theta}_{n,j}e_{1-j}$ \ \ \
$e_{1-j}\theta_n=\theta_{n,j}e_{1-j}$ \ \ \  $e_{1-j}(g_n-\chi_{cyc}(g_\infty))=v_{n,j}e_{1-j}$}
\]
so that $\tilde{\theta}_{n,j}\in \bbZ_p[\Gamma(n)]$ and
$\tilde{\theta}_{n,j}= v_{n,j}\theta_{n,j}$. Since $j\neq 0$, the augmentation of
$v_{n,j}=\omega(g_0)(\omega^{-j}(g_0)\gamma(n)-\kappa)$ lies in
$\bbZ_p^\times$ so that $v_{n,j}\in \bbZ_p[\Gamma(n)]^\times$ and
$\theta_{n,j}\in \bbZ_p[\Gamma(n)]$. Thus
$e_{1-j}\tilde{\theta}_\infty=v_{\infty,j}\theta_{\infty,j}e_{1-j}$ where
${\displaystyle v_{\infty,j}:=\lim_{\leftarrow} v_{n,j}\in \Lambda_\Gamma^\times}$ and
${\displaystyle \theta_{\infty,j}:=\lim_{\leftarrow} \theta_{n,j}\in \Lambda_\Gamma}$.
It follows that
$\Lambda_G^{(j)}/(\iota_\infty(\tilde{\theta}_\infty))^{(j)}=\Lambda_\Gamma e_j/\iota_\infty(\theta_{\infty,j})\Lambda_\Gamma e_j$.
Next, we identify $\Lambda_\Gamma$ as usual with $\Lambda:=\bbZ_p[[T]]$ by sending $\gamma$ to $1+T$. Then
$\theta_{\infty,j}$ goes to the unique power series  $f_j(T)\in\Lambda$ such that
\beq\label{eq: p-adic L function}
L_p(s,\omega^j\psi)=
f_j(\psi(\gamma(n))\inv\kappa^s-1)
\eeq
for all $s\in \bbZ_p$ and any character $\psi:\Gamma(n)\rightarrow\barbbQp^\times$ for any $n\geq 0$, where $L_p(s,\omega^j\psi)$ denotes the
$p$-adic $L$-function.
(See~\cite{Wash} pp.~119 and 122-123: take $\gamma(n)$ to be `$\gamma_n(1+q_0)$' for all $n$ so that
$\psi(\gamma(n))\inv$ equals `$\zeta_\psi$'  and check that
$\theta_{n,j}$ equals `$\xi_n(\omega^j)$' by~(\ref{eq: def of theta}).)
Thus $\Lambda_G^{(j)}/(\iota_\infty(\tilde{\theta}_\infty))^{(j)}\cong \Lambda/(f_j(\kappa(1+T)\inv-1))$.
But the `Main Conjecture' states in this case that $f_j(T)$ also equals $\charL(X_\infty^{(1-j)})$
(the characteristic power series of $X_\infty^{(1-j)}$ as a f.g.\ torsion $\Lambda$-module, defined up to a unit of
$\Lambda$). This is clearly equivalent to the two middle terms of~(\ref{eq: 4-term exact seq in case K=Q with j}) having the same
characteristic power series (up to a unit).
The Main Conjecture is, of course, proven (see \eg~\cite[Thm. 15.14]{Wash}, where characteristic {\em polynomials} are
used, for unicity). From
the multiplicativity of characteristic power series in exact sequences we deduce
\begin{thm}\label{prop: consequence of Main Conj} If $K=\bbQ$  then
\beq\label{eq: char p.s. and D}
\charL\left((A_\infty^{(j)})^\vee\right)=\charL(\Lambda_G^{(j)}/\fD_\infty^{(j)})\ \ \ \mbox{(up to a unit of $\Lambda$).}
\eeq
for all $j$ is even, $0\leq j\leq p-3$. \ePf
\end{thm}
Note that both sides of~(\ref{eq: char p.s. and D}) are units for $j=0$ and Greenberg's Conjecture is equivalent to the same for all even $j$. (See also
Prop.~\ref{prop: vandi and green}\ref{part: vandi and green (i)} for the R.H.S.)\vertsp\\
\rem\ Let $\fD_{\infty,j}$ be the ideal of $\Lambda_\Gamma$ determined
by $\fD_\infty^{(j)}=\fD_{\infty,j}e_j$. Identifying $\Lambda_\Gamma$ with $\Lambda$ (a noetherian U.F.D.)
we find that $\charL(\Lambda_G^{(j)}/\fD_\infty^{(j)})$
is simply the l.c.m.\ of any set of $\Lambda$-generators of $\fD_{\infty,j}$. Equation~(\ref{eq: char p.s. and D})
is equivalent to the statement that $\fD_{\infty,j}$ is contained with finite index in the principal ideal generated by
$\charL((A_\infty^{(j)})^\vee)$
and the Main Conjecture
for $K=\bbQ$ would follow by converse arguments if we could give an independent proof of this statement for all even $j$.
Consider the extra hypothesis that $(A_\infty^{(j)})^\vee$ is pseudo-isomorphic to a {\em cyclic} $\Lambda$-module $\Lambda/(c_j)$, so
$c_j=\charL((A_\infty^{(j)})^\vee)$. The inclusion $\fD_{\infty,j}\subset(c_j)$ then follows
from Theorem~\ref{thm: Stick in + part over Zp v2}
and the finiteness of the index should follow from~Cor.~\ref{prop: K=Q, Fittings, duals etc.}.
Without this hypothesis however,  a new ingredient would probably be required to re-prove the Main Conjecture by this route, perhaps an
`Euler Systems'-type elaboration of Theorem~\ref{thm: Stick in plus part mod pnpo}.
\vertsp\\
\rem\label{rem: annihilation by fS vs Gras-Oriat} By Thm.~\ref{prop: generation of fS} and the foregoing calculations, we have
$\iota_\infty(\fS_\infty)e_j=\iota_\infty(\theta_{\infty,j})\Lambda_\Gamma e_j$ for each $j\neq 0$,
so $\iota_\infty(\theta_{\infty,j})$ annihilates $X_\infty^{(j)}$ by Cor.~\ref{cor: containments and annihilation}~\ref{part: annihilation}.
For $n\geq 0$, let $d_{n,j}$ denote the image of $\iota_\infty(\theta_{\infty,j})$ in $\bbZ_p[\Gamma(n)]$. Since $\theta_{\infty,j}$
corresponds to $f_j(T)$ we find $\psi(d_{n,j})=f_j(\kappa\psi(\gamma(n))\inv-1)=L_p(1,\omega^j\psi)$ 
for every character $\psi:\Gamma(n)\rightarrow\barbbQp^\times$, giving the formula
$
d_{n,j}=\sum_{\psi}L_p(1,\omega^j\psi)e_\psi
$
(where $\psi$ ranges over all such characters and $e_\psi$ is the corresponding idempotent in $\barbbQp[\Gamma(n)]$).
Clearly, $d_{n,j}$ annihilates $(X_\infty^{(j)})_{\Gamma(n)}$ which is isomorphic to
$X_n^{(j)}\cong A_n^{(j)}$ (see Remark~\ref{rem: Gamma_n covts}, using $K=\bbQ$). This is a weakening of
the annihilation results of Gras and Oriat (see~\cite{Oriat}). The latter hold for more general real abelian fields
and even in the present case amount (more or less)
to the annihilation by $d_{n,j}$ of the much {\em bigger module}\ $\fX_n^{(j)}$. (Indeed, if $p$ divides the numerator of the
$j$th Bernoulli number, one can show
that $|\fX_n^{(j)}|$ is finite but unbounded as $n\rightarrow\infty$, whereas $A_n^{(j)}=\{0\}$ in all known cases.)
This, of course, correponds to the fact that $\iota_\infty(\fS_\infty)$ is usually
a much {\em smaller ideal}\ than $\fD_\infty$. Considerable generalisations of Gras' and
Oriat's annihilation results appear in~\cite{Burns-Barrett}. See also~\cite{Belliard-NQD}.
\vertsp\\
\rem\
There is a more familiar exact sequence
featuring both in a formulation of the Main Conjecture essentially due to Iwasawa in~\cite{Iwasawa on some mods...}
and in its proof by Rubin (see \eg~\cite[\S\S 15.4-7]{Wash}). Still in the case $K=\bbQ$ it reads (for each even $j$):
\beq\label{eq: more familiar exact seq}
0\rightarrow
(\hat{E}^1_\infty/\hat{C}^1_\infty)^{(j)}\longrightarrow
(U^1_\infty/\hat{C}^1_\infty)^{(j)}
\longrightarrow
\fX_\infty^{(j)}
\longrightarrow
X_\infty^{(j)}
\rightarrow 0
\eeq
Recall that $\hat{E}^1_n$ ({\em resp.}\ $\hat{C}^1_n$) denotes the {\em closure} in
$U^1(K_n)^+$ of the group $E^1_n$ ({\em resp.}\ $C^1_n$) which in turn consists of the embeddings
of those elements of $E_n$ ({\em resp.} of $C_n$, see Section~\ref{sec: d and D for K=Q}) which are congruent to $1$ modulo the unique prime above $p$
in $K_n^+$. Then  $\hat{E}^1_\infty$ ({\em resp.}\ $\hat{C}^1_\infty$) is
obtained by taking the projective limit w.r.t.\ norms. The middle map comes from
the map $\psi_\infty$ used in Section~\ref{sec: inertia subgps} (but here in the {\em plus} part). We assume for simplicity that $j\neq 0$
and compare the non-zero terms of this sequence with those of~(\ref{eq: 4-term exact seq in case K=Q with j}).
First, Iwasawa proved that $(U^1_\infty/\hat{C}^1_\infty)^{(j)}$ is $\Lambda$-isomorphic to
$\Lambda/(f_j(\kappa(1+T)\inv-1))$. 
Hence
\[
(U^1_\infty/\hat{C}^1_\infty)^{(j)}\cong \Lambda_G^{(j)}/(\iota_\infty(\tilde{\theta}_\infty))^{(j)}
\]
For the remaining terms we use
the notion of the {\em adjoint}\ $\alpha(M)$ of a torsion $\Lambda$-module $M$ (see~\cite[\S15.5]{Wash}).
In the special case where $M_{\Gamma_n}$ is {\em finite}\ for all $n\geq 0$ we have
$\alpha(M)\cong {\displaystyle \lim_\leftarrow}(M_{\Gamma_n})^\vee$, where the map
$(M_{\Gamma_{n+1}})^\vee\rightarrow (M_{\Gamma_n})^\vee$ is dual to the map
$M_{\Gamma_{n}}\rightarrow M_{\Gamma_{n+1}}$ given by multiplication by $1+\gamma_n+\gamma_n^2+\ldots+\gamma_n^{p-1}$.
(See~\cite{Iwasawa on Z_l extensions...}.) From the isomorphism  $(X_\infty)_{\Gamma_n}\rightarrow X_n$ we deduce
$\alpha(X_\infty^{(i)})\cong A_\infty^{(i),\vee}\ \forall\, i$. Also, $\alpha$ commutes with $\dag$ and~(\ref{eq: second for  Lambda_Gamma Structure})
gives $A_\infty^{(1-j),\vee}\cong (\fX_\infty^{(j)})^\dag$ for $j$ even. Therefore
\[
\alpha(X_\infty^{(j)})\cong A_\infty^{(j),\vee}\ \ \ \mbox{and}\ \ \ \alpha((X_\infty^{(1-j)})^\dag)\cong \fX_\infty^{(j)}.
\]
Finally, from $(\Lambda_G^{(j)}/\fD_\infty^{(j)})_{\Gamma_n}\cong(\bbZ_p[G_n]/\fD_n)^{(j)}$ and Thm.~\ref{prop: gp ring mod fD iso to E/C dual}
one deduces that $\alpha(\Lambda_G^{(j)}/\fD_\infty^{(j)})$ is the projective limit of the groups
$\Zpten(\tE_n/\tC_n)\cong (\Zpten E^1_n)/(\Zpten C^1_n)$ w.r.t.\
the norm maps. Now, since Leopoldt's Conjecture holds for $K_n$, there is an isomorphism $\Zpten E^1_n\rightarrow \hat{E}^1_n$ taking
$\Zpten C^1_n$ to $\hat{C}^1_n$. Hence $(\Zpten E^1_n)/(\Zpten C^1_n)\cong \hat{E}^1_n/\hat{C}^1_n$ and the compactness of $\hat{E}^1_n$,
gives
\[
\alpha(\Lambda_G^{(j)}/\fD_\infty^{(j)})\cong(\hat{E}^1_\infty/\hat{C}^1_\infty)^{(j)}
\]
Despite these relations between the terms of~(\ref{eq: more familiar exact seq}) and~(\ref{eq: 4-term exact seq in case K=Q with j}), it is
not obvious to the author that one sequence
follows directly from the other or even whether such neat relations are to be expected between the terms of
appropriately generalised sequences for any abelian $K$.

\end{document}